\documentclass[12pt]{amsart}
\usepackage{amsfonts}
\usepackage{latexsym,amsmath,amsthm,amssymb,amsfonts}

\numberwithin{equation}{section}
\allowdisplaybreaks

\def\endproof{$\hfill\Box$\\}
\def\s{\,\,\,\,}

\numberwithin{equation}{section}
\newtheorem{theorem}{Theorem}[section]
\newtheorem{lem}[theorem]{Lemma}
\newtheorem{thm}[theorem]{Theorem}

\newtheorem{defi}[theorem]{Definition}
\newtheorem{rem}[theorem]{Remark}

\def\dint{\displaystyle{\int}}

\newcounter{Cnumber}

\title[ ]
{Global Weak Solutions to Landau-Lifshitz Equations into Compact Lie Algebras}
\author[ ]
{Zonglin Jia,\quad\quad  Youde Wang}

\thanks{The authors are supported by NSFC grant (No.11731001).}

\date{}
\begin{document}

\maketitle

\begin{abstract}
In this paper, we consider a parabolic system from a bounded domain in a Euclidean space or a closed Riemannian manifold into a unit sphere in a compact Lie algebra $\mathfrak{g}$, which can be viewed as the extension of Landau-Lifshtiz (LL) equation and was proposed by V. Arnold. We follow the ideas taken from the work by the second author to show the existence of global weak solutions to the Cauchy problems of such Landau-Lifshtiz equations from an $n$-dimensional closed Riemannian manifold $\mathbb{T}$ or a bounded domain in $\mathbb{R}^n$ into a unit sphere $S_\mathfrak{g}(1)$ in $\mathfrak{g}$. In particular, we consider the Hamiltonian system associated with the nonlocal energy--{\it micromagnetic energy} defined on a bounded domain of $\mathbb{R}^3$ and show the initial-boundary value problem to such LL equation without damping terms admits a global weak solution. The key ingredient of this article consists of the choices of test functions and approximate equations.

\medskip
\break
{\textbf{Key words}}: Lie Algebra, Test Functions

\medskip
{\textbf{AMS Subject Classification}}: 35D30, 35G20, 35G25

\end{abstract}

\section{Introduction}
Let $\Omega$ be a bounded domain in the Euclidean space $\mathbb{R}^3$. In physics, the Landau-Lifshitz equation
$$u_t+u\times(\Delta_{\mathbb{R}^3} u-J(u))=0$$
was introduced by Landau and Lifshitz \cite{LL} as a model for the magnetization $u:\Omega\to S^2$ in a ferromagnetic material. The matrix $J:=\mbox{diag}(J_1, J_2, J_3)$ gives account of the anisotropy of the material. The equation describes the Hamiltonian dynamics corresponding to the Landau-Lifshitz energy
\[
E(u)=\frac{1}{2}\int_{\Omega}(|\nabla_{\mathbb{R}^3} u|^2+\lambda_1u_1^2+\lambda_3u_3^2)\,dx.
\]
where the $\nabla_{\mathbb{R}^3}$ denotes the gradient operator on $\mathbb{R}^3$ and $dx$ is the volume element. The two values of the characteristic numbers $\lambda_1:=J_2-J_1$ and $\lambda_3:=J_2-J_3$ are non-zero for biaxial ferromagnets, while $\lambda_1$ is chosen to be equal to $0$ in the case of uniaxial ferromagnets. The material is isotropic when $\lambda_1=\lambda_3=0$, and the Landau-Lifshitz equation reduces to the well-known Schr\"odinger map flow equation (\cite{DW, BIKT}).

The Landau-Lifshitz equation with dissipation, which can be written as
\[
u_t=-u\times\Delta_{\mathbb{R}^3} u +\alpha u\times u_t,\]
was proposed by Gilbert in 1955 \cite{G}. Here $\alpha>0$ is the damping parameter, which is characteristic of the material, and $\alpha$ is usually called the Gilbert damping coefficient. Hence the Landau-Lifshitz equation with damping term is also called the Landau-Lifshitz-Gilbert (LLG) equation in the literature.

Generally, in physics the Landau-Lifshitz functional is defined by
$$\mathcal{E}(u):=\int_{\Omega}\Phi(u)\,dx+\frac{1}{2}\int_{\Omega}|\nabla u|^2\,dx-\frac{1}{2}\int_{\Omega}h_d\cdot u\,dx.$$
In the above functional, the first and second terms are the anisotropy and exchange energies, respectively. $\Phi(u)\equiv \sum_{i=1}^3\lambda_iu_i^2$ is a real function on $S^2$, where $\lambda_i$ are real number. Besides, one also considers uniaxial materials with easy axis parallel to the OX-axis, for which $\Phi(u)= u_2^2 + u_3^2$. The last term is the self-induced
energy, and $h_d = -\nabla w$ is the demagnetizing field. The magnetostatic potential, $w$, solves the differential equation (the stray field equation)
$$\Delta w = \mbox{div}(u\chi_{\Omega})$$
in $\mathbb{R}^3$ in the sense of distributions, where $\chi_{\Omega}$ is the characteristic function of $\Omega$. The solution to this equation is
$$w(x) = \int_{\Omega} \nabla N(x-y)u(y)dy,$$
where $N(x) = -\frac{1}{4\pi|x|}$ is the Newtonian potential in $\mathbb{R}^3$. $\mathcal{E}(u)$ is a nonlocal energy as $h_d$ does not vanish.

In the absence of spin currents, the relaxation process of the magnetization distribution is described by the following
\begin{equation}\label{LLG}
u_t - \alpha u\times u_t=-u\times h,
\end{equation}
with $|u|=1$ and Neumann boundary condition:
$$\frac{\partial u}{\partial\nu}=0\s\s\s \mbox{on}\s \partial\Omega,$$
where $\nu$ represents the outward unit normal on $\partial\Omega$. The local field $h$ of $\mathcal{E}(u)$ is just
$$h:=-\frac{\delta\mathcal{E}(u)}{\delta u}=-\nabla_u\Phi + \Delta u + h_d.$$

In this paper, our first aim is to consider the following Landau-Lifshitz (LL) equation:
\begin{equation}\label{spin:0}
\left\{ \begin{aligned}
&\partial_tu =-u\times(\Delta u+h_d-\nabla_u\Phi),   \s\s   (x, t)\in\Omega\times\mathbb{R}^+,\\
& u(\cdot,0)=u_0:\Omega\longrightarrow S^2, \s\s\s\s\s\s\s \frac{\partial u}{\partial \nu}\big|_{\partial\Omega}=0.
\end{aligned} \right.
\end{equation}
In particular, the above system can be regarded as the Schr\"odinger flow corresponding to $\mathcal{E}(u)$ in the sense of the definition given in \cite{DW0}.
\medskip

Recall that Arnold and Khesin had ever proposed in \cite{AK} considering the so-called Landau-Lifshitz model associated with a Lie algebra(see p333-335 of \cite{AK}). For more details, we refer to \cite{AK} and Section 2 of this paper. In fact, Ding, Wang and Wang in  \cite{DWW} have ever discussed the existence of global weak solution to the Landau-Lifshitz systems from a closed Riemannian manifold into the unit sphere of a compact Lie algebra.

In this paper, we are intend to extend these models to the case the unknown functions $u$ are of Lie algebra value. Let $(\mathfrak{g}, [\cdot, \cdot])$ denote an $n$-dimensional compact semisimple Lie algebra associated with a compact semisimple Lie group $G$ and $S_\mathfrak{g}(1)$ denote the unit sphere in $\mathfrak{g}$ centered at the origin. Now we are going to generalize the above model to the case the target manifold is $S_\mathfrak{g}(1)$. First, we would like to consider the following Landau-Lifshitz-Gilbert system with Lie algebra value:
\begin{equation}\label{spin:1}
\left\{ \begin{aligned}
&\partial_tu-\alpha[u,\,\partial_tu]=-[u,\,\Delta u+h_d-\nabla_u\Phi],              & (x, t)\in\Omega\times\mathbb{R}^+,\\
& u(\cdot,0)=u_0:\Omega\longrightarrow S_{\mathfrak{g}}(1), \s\s \frac{\partial u}{\partial \nu}\big|_{\partial\Omega}=0.
\end{aligned} \right.
\end{equation}
Here, $\Omega$ is a bounded domain in Euclidean space $\mathbb{R}^n$ with boundary $\partial\Omega$, $\nu$ is the outward unit normal of $\partial\Omega$ and $h_d=-\nabla w$ where $w$ satisfies $$\Delta w= \mbox{div}(u\chi_{\Omega}).$$
It is worthy to point out that $\dim(\Omega)=\dim(\mathfrak{g})$ is needed. Otherwise, the definition of $h_d$ does not make sense.

\medskip
On the other hand, if the Landau-Lifshitz energy is defined by
\[
E(u)=\frac{1}{2}\int_{\Omega}f\cdot|\nabla_{\mathbb{R}^3} u|^2\,dx,
\]
where $f$ is a positive smooth function which is usually called coupling function. Then, the corresponding Landau-Lifshitz equation is
\[
u_t - \alpha u\times u_t=-u\times(f\triangle_{\mathbb{R}^3} u+\nabla_{\mathbb{R}^3} f\cdot\nabla_{\mathbb{R}^3} u).
\]
This system is usually called inhomogeneous Landau-Lifshitz equation which was discussed by some physicists and mathematicians(see \cite{B, DPL}). We can also make the following extension of the inhomogeneous Landau-Lifshitz equation. Let $(\mathbb{T}, \mathfrak{h})$ be an $n$-dimensional closed Riemannian manifold equipped with a metric $\mathfrak{h}=(h_{ij})$.
\begin{equation}\label{spin:2'}
\left\{ \begin{aligned}
&\partial_tu-\alpha[u,\,\partial_tu]=-[u,\,\Delta_f u-\nabla_u\Phi],              & x\in\mathbb{T},\\
& u(\cdot,0)=u_0:\mathbb{T}\longrightarrow S_{\mathfrak{g}}(1), &
\end{aligned} \right.
\end{equation}
where $u:\mathbb{T}\times\mathbb{R}^+\rightarrow S_{\mathfrak{g}}(1)$ is an unknown mapping and
$$\Delta_fu\equiv f\Delta u+\nabla f\cdot\nabla u.$$
Here, $\Delta$ and $\nabla$ denote respectively Laplace-Beltrami operator and gradient operator on $\mathbb{T}$ with respect to the metric $\mathfrak{h}$.

More generally, from the viewpoint of mathematics we may consider the following equation:
\begin{equation}\label{lie:1}
\left\{
\begin{aligned}
&\alpha_0u_t+\alpha[u,\,u_t]=[u,\,f\Delta u+\nabla f\cdot\nabla u]+F(x,t,u),\\
&u(\cdot ,0)=u_0:\mathbb{T}\rightarrow S_\mathfrak{g}(1),\s\s\s u_0\in H^{1}(\mathbb{T},S_\mathfrak{g}(1)).
\end{aligned}
\right.
\end{equation}
Here $\alpha_0>0$ and $\alpha\geqslant0$ are two constants. $f:\mathbb{T}\rightarrow \mathbb{R}^+$ is a $C^1$-function and $F$ is a $C^1$-smooth mapping from $\mathbb{T}\times\mathbb{R}^+\times\mathfrak{g}$ to $\mathfrak{g}$ which satisfies the following
\[
\langle F(x,t,z),z\rangle=0,\,\,\,\s \forall\,\,z\in\mathfrak{g}.
\]
We call the above as generalized inhomogeneous Landau-Lifshitz equation or GILL equation for short.

\medskip
In recent years, there has been lots of interesting studies for the Landau-Lifshitz equation, concerning its existence, uniqueness and regularities of various kinds of solutions. Before moving on to the next step, we list only a few of the literature that are closely related to our work in the present paper.

First, let us recall some results on Landau-Lifshitz equation. In the case $\alpha>0$, the existence and non-uniqueness of weak solutions to the LLG equation goes back to \cite{AS, V}. For $\Omega$ is a bounded domain in $\mathbb{R}^3$, Carbou and Fabrie studied a model of ferromagnetic material governed by a nonlinear Landau-Lifschitz equation coupled with Maxwell equations and proved the local existence of a unique strong solutions in \cite{CF} (also see \cite{CR, DLW, DG}). In fact, in two space dimensions and for sufficiently small initial data, the strong solution is global in time \cite{CF}. For general initial data, the two-dimensional solution may develop finitely many point singularities after finite time; see \cite{Hap} for a discussion. Later, Tilioua \cite{T} employed the penalized method to show the existence of the weak solution to LLG with spin current.

In case $\alpha=0$, $f=1$, and $[\cdot\, ,\,\cdot]$ is just the cross product in $\mathbb{R}^3$, the existence and uniqueness of smooth solutions to LL equation goes back to \cite{SSB, ZGT}. For the system (\ref{lie:1}) Wang has established the existence of global solution to LL equation without Gilbert damping term defined on a closed Riemannian manifold in \cite{W}. Later, in \cite{DWW} the authors proved the existence of global weak solutions to (\ref{lie:1}) on a closed Riemannian manifold $\mathbb{T}$ or $\mathbb{T}=\mathbb{R}^n$ for the case $\alpha=0$, $f=1$ and $[\cdot,\,\cdot]$ is the Lie bracket of a compact Lie algebra.

The essential difference between $(\mathbb{R}^3, \times)$ and $(\mathfrak{g}, [\cdot, \cdot])$ lies on that, for $|u|=1$, in $(\mathbb{R}^3, \times)$ $u\times(u\times v)=-v$ for any $v\in T_u S^2$, where $S^2$ is the unit sphere of $\mathbb{R}^3$ and $T_u S^2$ is its tangent space at $u$, while $[u, [u, v]]=-v$ is not true generally for any $v\in T_u\left(S_{\mathfrak{g}}(1)\right)$ in $(\mathfrak{g}, [\cdot, \cdot])$. Therefore, $u_t - \alpha u\times u_t=-u\times h$ is equivalent to $u\times u_t + \alpha u_t= h$ in $(\mathbb{R}^3, \times)$, but generally $$u_t - \alpha [u, u_t]=-[u, h]$$ is not equivalent to $$[u, u_t] - \alpha[u, [u, u_t]]= -[u, [u, h]]$$ in $(\mathfrak{g}, [\cdot, \cdot])$. Hence, for LL equation with Lie algebra value the well-known Ginzburg-Landau penalized method as in \cite{CF, CR} is not again effective. Even for the heat flows of harmonic map  $$u_t=\tau(u)=\Delta u +|\nabla u|^2u,\s\s\s u: \mathbb{T}\times\mathbb{R}^+\to S_{\mathfrak{g}}(1)\subset\mathfrak{g}$$
such a penalized method adopted by Chen in \cite{Ch} is not again valid.

While one of the authors in \cite{W} has ever consider a different method to approach the existence of weak solution to Schr\"odinger flow for maps from a closed manifold or a bounded domain in Euclidean space. It seems that the method in \cite{W} is more effective for the present situation than the penalized method mentioned in the above.

In this paper, we follow the ideas in \cite{W} to approach the existence problems of systems (\ref{spin:1}) and (\ref{lie:1}). In particular, we focus on the case $\alpha=0$, i.e. Landau-Lifshitz equation without dissipation which associates with the nonlocal energy--{\it micromagnetic energy}. In the present situation, the nonlocal property of energy results in some new difficulties. By the elliptic regularity theory and delicate analysis we can overcome these obstructions and still obtain some uniform a priori estimates with respect to small enough $\alpha>0$. In the forthcoming paper \cite{JW1}, we will also consider the existence for the global weak solutions to Landau-Lifshitz systems with spin-polarized transport. Before stating our main results, we need to elucidate some definitions on the weak solution to the above equations which are given in Section 2. Now we present our main results.

\begin{thm}\label{thm}
Let $\Omega$ be  a bounded domain of an $n$-dimensional Euclidean space and $\mathfrak{g}$ be a $\dim{\Omega}$-dimensional compact Lie algebra. Suppose that $u_0$ belongs to $H^{1}(\Omega,\mathfrak{g})$ and $|u_0|=1$ a.e. $\Omega$. Then (\ref{spin:1}) admits a global weak solution with initial value $u_0$, provided $\Phi$ is $C^2$-smooth. In the case $\alpha=0$, $u\in L^\infty_{loc}(\mathbb{R}^+, H^1(\Omega, S_{\mathfrak{g}}(1)))$. In the case $\alpha >0$, $u\in W^{1,1}_2(\Omega\times[0, T],S_{\mathfrak{g}}(1))$ for any $T>0$.
\end{thm}

It is well-known that $(\mathbb{R}^3, \times)$ is a compact Lie algebra. As the first direct corollary of the above theorem, for LL equation or Schr\"odinger flow on $S^2$ we have
\begin{thm}\label{cor1}
Let $\Omega$ be a bounded domain of $\mathbb{R}^3$. Suppose that $u_0$ belongs to $H^{1}(\Omega,\mathbb{R}^3)$ with $|u_0|=1$ a.e. on $\Omega$. Then the initial boundary problem (\ref{spin:0}) of LL equation admits a global weak solution belonging to $L^\infty_{loc}(\mathbb{R}^+, H^1(\Omega, S^2))$ with initial value $u_0$, provided $\Phi$ is $C^2$-smooth.
\end{thm}

As the second direct corollary of the above theorem (\ref{thm}) which has been established essentially in \cite{GW} by a different method, we also have
\begin{thm}\label{cor2}
Let $\Omega$ be a bounded domain of $\mathbb{R}^3$. Suppose that $u_0$ belongs to $H^{1}(\Omega,\mathbb{R}^3)$ with $|u_0|=1$  a.e. on $\Omega$. Then, for any $T>0$ the initial boundary problem (\ref{LLG}) of LLG admits a global weak solution belonging to $W^{1,1}_2(\Omega\times[0, T],S^2)$ with initial value $u_0$, provided $\alpha>0$ and $\Phi$ is $C^2$-smooth.
\end{thm}

For the generalized inhomogeneous Landau-Lifshitz equation or (GILL for short) we obtain the following

\begin{thm}\label{thm1}
Let $(\mathbb{T}, \mathfrak{h})$ be an $n$-dimensional closed manifolds equipped with a metric $h$ and $\mathfrak{g}$ be a $m$-dimensional compact Lie algebra. Assume that  $\alpha \geq 0$, $F(x, t, z): \mathbb{T}\times\mathbb{R}^+\times S_{\mathfrak{g}}(1)\to \mathfrak{g}$ is $C^1$-smooth and $f\in C^1(\mathbb{T})$ with $\min_{x\in\mathbb{T}}f(x)>0$. Then (\ref{lie:1}) admits a global weak solution with initial value $u_0$ provided $u_0$ belongs to $H^{1}(\Omega, S_{\mathfrak{g}}(1))$. More precisely, in the case $\alpha_0>0$ and $\alpha>0$ the weak solution $u\in W^{1,1}_2(\mathbb{T}\times[0, T], S_{\mathfrak{g}}(1))$ for any $T>0$; in the case $\alpha_0>0$ and $\alpha=0$ the weak solution $u\in L^\infty(\mathbb{R}^+,H^1(\mathbb{T}, S_{\mathfrak{g}}(1)))$.
\end{thm}

The remainder of the article is organized as follows: In section 2 we recall some fundamental notions and summarize some known facts which will be used in this paper. In Section \ref{section3}, we will provide a proof of Theorem (\ref{thm}). An auxiliary approximate equation to (\ref{spin:1}) is chosen and an approximate solution to the auxiliary approximate equation is constructed based on a Galerkin approximation, and the necessary a priori estimates in order to guarantee the desired convergence are obtained. In Section \ref{section2}, by almost the same way as in Section \ref{section3} we consider the global well-posedness of generalized inhomogeneous LL equation and aim at showing Theorem \ref{thm1}.

\section{Preliminary}
First, we summarize some fundamental facts on compact Lie algebra. One may define a compact Lie algebra either as the Lie algebra of a compact Lie group, or as a real Lie algebra whose Killing form is negative definite. In this paper we always assume that $\mathfrak{g}$ is a Lie algebra whose Killing form is negative definite. It is well-known that, if the Killing form of a Lie algebra is negative definite, then the Lie algebra is the Lie algebra of a compact semisimple Lie group $G$. In general, the Lie algebra of a compact Lie group decomposes as the Lie algebra direct sum of a commutative summand (for which the corresponding subgroup is a torus) and a summand on which the Killing form is negative definite. It is well-known that there always is an $\mbox{Ad}(G)$-invariant inner product induced by the nondegenerate Killing form on $\mathfrak{g}$, denoted by $\langle\cdot,\cdot\rangle$(sometimes we omit it and sometimes we denote it by ``$\cdot$"), such that for any $X,Y,Z\in\mathfrak{g}$
\[
\langle Y,[X,Z]\rangle+\langle[X,Y],Z\rangle=0,
\]
where $[\cdot\,,\cdot]$ is the Lie bracket. So, it follows that there holds true
\[
\langle X,[X,Z]\rangle=0.
\]
In fact, $\mathbb{R}^3$ with cross product is just a compact Lie algebra corresponding to $SO(3)$. For the details we refer to the chapter 4 of \cite{H} and \cite{XWY,DWW}.

For instance, let $G$ be a finite-dimensional matrix group with a nondegenerate Killing form $\langle A, B\rangle=-\mbox{tr}(AB)$ for $A, B\in G$; i.e., a reductive group(one can think of $SO(3)$ or the group of all nondegenerate matrices $GL(n)$ and let $\mathfrak{g}$ be the corresponding Lie algebra). One defines the loop group $\tilde{G}$ is the group of $G$-value functions on the circle $\tilde{G}=C^\infty(S^1, G)$ with pointwise multiplication. The corresponding loop Lie algebra $\tilde{\mathfrak{g}}$ is the Lie algebra of $\mathfrak{g}$-value functions on the circle with pointwise commutator.

One had shown that the Landau-Lifshitz equation associated to a Lie algebra $\mathfrak{g}$ is the Euler equation corresponding to the loop group $\tilde{G}$ with the quadratic Hamiltonnian functional
$$E(u)=-\frac{1}{2}\int_{S^1}\mbox{tr}(\partial_x u)^2\,dx$$
on the dual space $\tilde{\mathfrak{g}}^*$, where $\partial_x u$ is $\mathfrak{g}$-value function, and ``$\mbox{tr}$" stands for the trace in the matrix algebra $\mathfrak{g}$ (Theorem 3.17, p334 in \cite{AK}).

The compact Lie algebras are classified and named according to the compact real forms of the complex semisimple Lie algebras. These are:

\noindent$\bullet$ ${A}_n: \mathfrak{su}_{ n + 1}$ corresponding to the special unitary group (properly, the compact form is $PSU$, the projective special unitary group);

\noindent$\bullet$ $B_n: \mathfrak {so}_{2n+1}$ corresponding to the special orthogonal group (or $\mathfrak {O}_{2n+1}$ corresponding to the orthogonal group);

\noindent$\bullet$ $C_n: \mathfrak {sp}_{n}$ corresponding to the compact symplectic group; sometimes written $\mathfrak {usp}_{n}$;

\noindent$\bullet$ $D_n: \mathfrak {so}_{2n}$ corresponding to the special orthogonal group (or $\mathfrak {O}_{2n}$ corresponding to the orthogonal group) (properly, the compact form is $PSO$, the projective special orthogonal group);

\noindent$\bullet$ Compact real forms of the exceptional Lie algebras $E_6$, $E_7$, $E_8$, $F_4$, $G_2$.

\medskip

Next, we recall some notions and notations on manifolds. Let $(\mathcal{M}, g)$ and $(\mathcal{N},\tilde{g})$ be two Riemannian manifolds and $\mathcal{N}$ is embedded isometrically in $\mathbb{R}^K$. The energy of a map $u$ from $\mathcal{M}$ into $\mathcal{N}$ is defined by
$$E(u)=\frac{1}{2}\int_\mathcal{M}|\nabla u|^2\,d\mathcal{M}.$$
The tension field of a map from $(\mathcal{M}, g)$ into $(\mathcal{N}, \tilde{g})$ is given by
$$\tau(u)=\Delta_\mathcal{M} u +A(u)(\nabla u, \nabla u),$$
where $\Delta_\mathcal{M}$ denotes the Laplace-Beltrami operator defined on $(\mathcal{M}, g)$ and $A(\cdot,\cdot)$ is the second fundamental form of $u$ in $(\mathcal{N},\tilde{g})$. In local charts, it can be written as
$$\tau(u)^\alpha=\Delta_\mathcal{M} u_\alpha + g^{ij}(x)\Gamma^\alpha_{\beta\gamma}(u)\frac{\partial u_\beta}{\partial x_i}\frac{\partial u_\gamma}{\partial x_j}.$$
Here, $\Gamma^\alpha_{\beta\gamma}$ is the Christoffel symbols of $(\mathcal{N}, \tilde{g})$.

In local coordinates $(x^1,...,x^n)$ on $(\mathcal{M}, g)$ and coordinates $(u_1,...,u_m)$ on $\mathfrak{g}$,
\[
\Delta_\mathcal{M} u=(\Delta_\mathcal{M} u_1, ..., \Delta_\mathcal{M} u_m)
\]
where
\[
\Delta_\mathcal{M} u_\alpha=\frac{1}{\sqrt{g}}\frac{\partial}{\partial x^i}\big(g^{ij}\sqrt{g}\frac{\partial u_\alpha}{\partial x^j} \big),\,\,\,\,\quad\mbox{for}\quad\alpha=1,2,...,m,
\]
and $(g^{ij})$ is the inverse of $(g_{ij})$. For convenience we always denote $\Delta_\mathcal{M}$ by $\Delta$.

On the other hand, for a smooth function $f$ defined on $\mathcal{M}$, we denote
\[
\nabla f\cdot\nabla u=(\nabla f\cdot\nabla u_1,\,...\, , \nabla f\cdot\nabla u_m)
\]
where
\[
\nabla f\cdot\nabla u_\alpha=\frac{\partial f}{\partial x^i}\frac{\partial u_\alpha}{\partial x^p}g^{ip},\,\,\,\,\quad\mbox{for}\quad \alpha=1,2,...,m.
\]

It is easy to see that the flow $u_t=[u, \, \Delta u]$ from $\mathcal{M}$ into a unit sphere of $\mathfrak{g}$ conserves $E(u)$ if the flow is smooth. Moreover, the the flow can also be written as
$$u_t=[u,\,\tau(u)]$$
since $\tau(u)=\Delta_\mathcal{M}u + |\nabla u|^2u$ in the present situation. It is just the Schr\"odinger flow if $\mathfrak{g}$ is replaced by $\mathbb{R}^3$.

We define the Sobolev spaces of the functions with compact Lie algebra value $$L^2(\mathcal{M}, \mathfrak{g})=\{u:\, \int_\mathcal{M}|u(x)|^2\,d\mathcal{M}<\infty\},$$ $$H^1(\mathcal{M}, \mathfrak{g})=\{u:\, |u|, |\nabla u|\in L^2(\mathcal{M})\}$$ and for maps from $\mathcal{M}$ into $\mathcal{N}$ by
$$H^1(\mathcal{M}, \mathcal{N})=\{u:\,u\in H^1(\mathcal{M}, \mathbb{R}^K),\, u(x)\in \mathcal{N} \s a.e.\,\, x\in\mathcal{M}\}.$$
Moreover, we define
$$W^{r,s}_p(\mathcal{M}, \mathcal{N})=\{u:\,  u\in W^{r,s}_p(\mathcal{M}, \mathbb{R}^K), \,u(x)\in \mathcal{N} \s a.e.\,\, x\in\mathcal{M}\}.$$
Similarly, we define $H^1(\mathbb{T}, S_{\mathfrak{g}}(1))= \{u:\, u\in H^1(\mathbb{T}, \mathfrak{g}), \,u(x)\in S_{\mathfrak{g}}(1) \s a.e.\,\, x\in\mathbb{T}\}$. And we say that $u\in L^\infty_{loc}(\mathbb{R}^+,H^1(\mathbb{T}, S_{\mathfrak{g}}(1)))$ means that, for any $T>0$, $u\in L^\infty([0, T],\,H^1(\mathbb{T}, S_{\mathfrak{g}}(1)))$.

\medskip
We will use the property of the following operator which is defined by
$$h_d(u)=-\nabla(\nabla N\ast u): L^2(\Omega,\mathbb{R}^n)\rightarrow L^2(\mathbb{R}^n,\,\mathbb{R}^n),$$
where $N(|x-y|)$ is the classical Newton potential and $n=\dim(\Omega)$. Then, for $u,\, \tilde{u}\in L^2(\Omega,\mathbb{R}^n)$ there hold
$$\int_{\mathbb{R}^n}|h_d(u)|^2\,dx\leq\int_{\Omega}|u|^2\,dx,$$
and
$$\int_{\mathbb{R}^n}|h_d(u)-h_d(\tilde{u})|^2\,dx\leq\int_{\Omega}|u-\tilde{u}|^2\,dx.$$
The fields $h_d(u)$ can be defined equivalently by $$h_d(u)=-\nabla w$$ where
$$\Delta w = \mbox{div}(u\chi_{\Omega})\s\s \mbox{in} \s \mathbb{R}^n$$
in the sense of distributions.

Multiplying the equation by any $v\in H^1(\mathbb{R}^n)$ and integrating by parts, we obtain
$$\int_{\mathbb{R}^n}\nabla w\cdot\nabla v=\int_{\Omega}u\cdot\nabla v.$$
Takeing $v=w$ in the above identity we get
$$\int_{\mathbb{R}^n}|\nabla w|^2=\int_{\Omega}u\cdot\nabla w\leq\left(\int_{\mathbb{R}^n}|\nabla w|^2\right)^{\frac{1}{2}}
\left(\int_{\Omega}|u|^2\right)^{\frac{1}{2}}.$$
It follows
$$\int_{\mathbb{R}^n}|h_d(u)|^2=\int_{\mathbb{R}^n}|\nabla w|^2\leq \int_{\Omega}|u|^2.$$

\medskip
In fact, the following lemma was shown in \cite{CF} and \cite{GW} although they only need to consider the case $\dim(M)=n=3$ therein. For more details we refer to page 196 in \cite{La}.

\begin{lem}\label{lem} For any $u, \, \tilde{u}\in L^2(\Omega,\mathbb{R}^n)$ with $n=\dim(\Omega)$, the operator $h_d$ satisfies
$$\int_{\mathbb{R}^n}|h_d(u)|^2\,dx\leq\int_{\Omega}|u|^2\,dx,$$
and
$$\int_{\mathbb{R}^n}|h_d(u)-h_d(\tilde{u})|^2\,dx\leq\int_{\Omega}|u-\tilde{u}|^2\,dx.$$
Moreover, if $u$ belongs to $W^{1,p}(\Omega)$ and $p\in (0, +\infty)$, the restriction of $h_d(u)$ to $\Omega$ belongs to $W^{1,p}(\Omega)$ and there exists a constant $C$ such that
$$\|h_d(u)\|_{W^{1,p}(\Omega)}\leq C\|u\|_{W^{1,p}(\Omega)}.$$
\end{lem}

Now we give the definitions of weak solutions:
\begin{defi}
In the case $\alpha >0$, we say $u$ is a global weak solution of equation $(\ref{spin:1})$ with initial data $u_0$ if

1. $u\in L^{\infty}(\mathbb{R}^+,H^1(\Omega,\mathfrak{g}))$, $\partial_tu\in L^2(\Omega\times\mathbb{R}^+,\mathfrak{g})$ and $|u|=1$ a.e. on $\mathbb{R}^+\times\Omega$,

2. for all $\varphi\in C^{\infty}(\Omega\times[0, T],\mathfrak{g})$, we have
\begin{eqnarray*}
&&\int_0^Tdt\int_{\Omega}\partial_tu\cdot\varphi\,d\Omega-\alpha\cdot\int_0^Tdt\int_{\Omega}[u,\partial_tu]\cdot\varphi\,d\Omega\\
&=&\sum\limits_{p=1}^n\int_0^Tdt\int_{\Omega}[u,\nabla_{p}u]\cdot\nabla_{p}\varphi\,d\Omega-\int_0^Tdt\int_{\Omega}[u,h_d-\nabla_u\Phi]\cdot\varphi\,d\Omega
\end{eqnarray*}
where $h_d:=-\nabla(\nabla N\ast u)$ and $N$ is the Newtonian potential in $\mathbb{R}^n$.

3. $u(0,x)=u_0(x)$ in the trace sense.

\medskip
In the case $\alpha =0$, we say that $u\in L^{\infty}(\mathbb{R}^+,H^1(\Omega,\mathfrak{g}))$ and $|u|=1$ a.e. on $\Omega\times\mathbb{R}^+$ is a global weak solution of equation $(\ref{spin:1})$ with initial data $u_0$ if, for all $\varphi\in C^{\infty}(\Omega\times[0,T],\mathfrak{g})$, we have
\begin{eqnarray*}
&&\int_{\Omega}\langle u(T),\varphi(T)\rangle\,d\Omega - \int_{\Omega}\langle u_0,\varphi(0)\rangle\,d\Omega + \int_0^Tdt\int_{\Omega}[u,h_d-\nabla_u\Phi]\cdot\varphi\,d\Omega\\
&=&\sum\limits_{p=1}^n\int_0^Tdt\int_{\Omega}[u,\nabla_{p}u]\cdot\nabla_{p}\varphi\,d\Omega + \int_{0}^{T}\int_{\Omega}\langle u,\varphi_t\rangle\,d\Omega dt
\end{eqnarray*}
where $h_d:=-\nabla(\nabla N\ast u)$. Moreover, $u(0,x)=u_0(x)$ in the trace sense.
\end{defi}

\begin{defi}
In the case $\alpha>0$, we say $u$ is the weak solution of (\ref{lie:1}) with initial value $u_0$, if it is a function belonging to $L^{\infty}([0,T],H^1(\mathbb{T},S_\mathfrak{g}(1)))\bigcap W^{1,1}_2([0,T]\times\mathbb{T},\mathfrak{g})$ which satisfies the relation:
\[
\aligned
&\alpha_0\int_{\mathbb{T}}\langle u(T),\varphi(T)\rangle\,d\mathbb{T}-\alpha_0\int_{\mathbb{T}}\langle u_0,\varphi(0)\rangle\,d\mathbb{T}+\alpha\int_{0}^{T}\int_{\mathbb{T}}\langle[u,u_t],\varphi\rangle\,d\mathbb{T}dt\\
=&\alpha_0\int_{0}^{T}\int_{\mathbb{T}}\langle u,\varphi_t\rangle\,d\mathbb{T}dt-\sum\limits_{p=1}^n\int_{0}^{T}\int_{\mathbb{T}}\langle[u,f\nabla_{e_p} u],\nabla_{e_p}\varphi\rangle\,d\mathbb{T}dt+\int_{0}^{T}\int_{\mathbb{T}}\langle F(x,t,u),\varphi\rangle\,d\mathbb{T}dt,
\endaligned
\]
for all $\varphi\in C^{\infty}(\mathbb{T}\times[0,T],\mathfrak{g})$. Here $\{e_p:1\leqslant p\leqslant n\}$ is a local orthonormal frame on $\mathbb{T}$.

\medskip
In the case $\alpha=0$, $u$ is a weak solution to (\ref{lie:1}) if $u\in L^{\infty}([0,T],H^1(\mathbb{T},S_\mathfrak{g}(1)))$ satisfies the following relation:
\[
\aligned
&\alpha_0\int_{\mathbb{T}}\langle u(T),\varphi(T)\rangle\,d\mathbb{T}-\alpha_0\int_{\mathbb{T}}\langle u_0,\varphi(0)\rangle\,d\mathbb{T}+\sum\limits_{p=1}^n\int_{0}^{T}\int_{\mathbb{T}}\langle[u,f\nabla_{e_p} u],\nabla_{e_p}\varphi\rangle\,d\mathbb{T}dt\\
=&\alpha_0\int_{0}^{T}\int_{\mathbb{T}}\langle u,\varphi_t\rangle\,d\mathbb{T}dt+\int_{0}^{T}\int_{\mathbb{T}}\langle F(x,t,u),\varphi\rangle\,d\mathbb{T}dt,
\endaligned
\]
for all $\varphi\in C^{\infty}(\mathbb{T}\times[0,T],\mathfrak{g})$.
\end{defi}

\section{Initial-Boundary Value Problems of LLG Equations}\label{section3}
In this section, we will show the existence of weak solutions to the following Landau-Lifshitz-Gilbert system with Lie algebra value:
\begin{equation*}
\left\{ \begin{aligned}
&\partial_tu-\alpha[u,\,\partial_tu]=-[u,\,\Delta u+h_d-\nabla_u\Phi],              & x\in\Omega,\\
& u(\cdot,0)=u_0:\Omega\longrightarrow S_{\mathfrak{g}}(1), \s\s \frac{\partial u}{\partial \nu}\big|_{\partial\Omega}=0.
\end{aligned} \right.
\end{equation*}
We will see that the proof of Theorem \ref{thm} is based on a combination of a Galerkin approximation for the LLG equation (\ref{spin:1}), delicate choices of approximation equation and the picking of test functions. For this end, the first step we need to choose the following approximation equation:
\begin{equation}\label{A0}
\left\{ \begin{aligned}
&\partial_t u  = \varepsilon\Delta u+\left[\mathfrak{J}(u) ,\,\alpha\partial_t u-\Delta u-h_d(u)+\nabla_u\Phi\left(\mathfrak{J}(u)\right)\right], & x\in\Omega,\\
& u(\cdot,0)=u_0:\Omega\longrightarrow S_{\mathfrak{g}}(1), \s\s\s \frac{\partial u}{\partial \nu}\big|_{\partial\Omega}=0,
\end{aligned} \right.
\end{equation}
where
$$\mathfrak{J}(u)=\frac{u}{\max\{|u|,1\}}.$$
It should be pointed out that the above $\Phi(u)$ has been extended to the closed ball $\overline{B}_{\mathfrak{g}}(1)\subset\mathfrak{g}$. In fact, we can extend $\Phi(z)$ by
\begin{equation*}
\tilde{\Phi}(z)=
\left\{
\begin{aligned}\zeta(|z|^2)\Phi\left(\frac{z}{\max\{\delta_0,|z|\}}\right),\s\s |z|^2>\delta_0,\\
0, \s\s |z|^2\leq\delta_0,
\end{aligned}\right.
\end{equation*}
where $\zeta(t):[0, 1]\to [0, 1]$ is a $C^2$-smooth function with $\zeta(t)\equiv 0$ on $[0, \, 2\delta_0]$ ($2\delta_0<1$) and $\zeta(1)=1$. It is easy to see that $\tilde{\Phi}$ is $C^2$-smooth on $\overline{B}_{\mathfrak{g}}(1)$. For simplicity, we still denote $\tilde{\Phi}$ by $\Phi$.

\subsection{Galerkin Approximation: A Priori Estimates.}

Let $\lambda_i(i=1,2,\cdots)$ ($0=\lambda_1\leq\lambda_2\leq\cdots\leq\lambda_i\leq\cdots$) be the eigenvalues of the operator $-\Delta$ on the domain
$$X:=\left\{\omega\in H^2(\Omega):\frac{\partial\omega}{\partial\nu}\Big|_{\partial\Omega}=0\right\},$$
and let $\{\omega^i: i=1, 2, \cdots\}$ be the orthonormal basis of the corresponding eigenfunctions. That is to say,
\begin{equation*}
\left\{ \begin{aligned}
         &-\Delta\omega^i = \lambda_i\cdot\omega^i, \\
                & \frac{\partial\omega^i}{\partial\nu}\Big|_{\partial\Omega}=0.
                          \end{aligned} \right.
\end{equation*}

According to Galerkin approximation, we consider the approximate solutions to the auxiliary equation (\ref{A0}) as follows
\begin{eqnarray*}
u^N(x,t):=\sum\limits_{i=1}^{N}\beta_i^N(t)\omega^i(x).
\end{eqnarray*}
Here $\{\beta_i^N(t)\}$ are unknown functions which take values in $\mathfrak{g}$ and assumed to satisfy the following ordinary differential equation:
\begin{equation}\label{1}
\left\{
\begin{array}{llll}
\aligned
\frac{d\beta_i^N}{dt}=&\alpha\sum\limits_{k=1}^N\int_{\Omega}\left[\mathfrak{J}(u^N),\frac{d\beta^N_k}{dt}
\right]\omega^k\omega^i\,dx +\varepsilon\int_{\Omega}\sum\limits_{k=1}^{N}\beta_k^N(-\lambda_k\omega^k)\omega^i\,dx\\
 &+\int_{\Omega}\left(\left[\mathfrak{J}(u^N),\sum\limits_{k=1}^{N}\beta_k^N\cdot(\lambda_k\omega^k)\right] +\left[\mathfrak{J}(u^N),(\nabla_u\Phi)\left(\mathfrak{J}(u^N)\right)\right]\right)
\omega^i\,dx\\
 &+\int_{\Omega}\left[\mathfrak{J}(u^N),\,\sum\limits_{k=1}^{N}\beta_k^N\cdot\nabla(\nabla N\ast\omega^k)\right]\omega^i\,dx,\\
 \beta_j^N(0)=&\int_{\Omega}u_0\cdot\omega^j\,dx.
\endaligned
\end{array}
\right.
\end{equation}
It is easy to see that (\ref{1}) can be written as
\[
(\mbox{Id} + A(\beta))\frac{d\beta}{dt}=B(\beta).
\]
Here $\beta=(\beta_1^N,\beta_2^N,...,\beta_N^N)^T$, $\mbox{Id}$ is the unit matrix, and $A(\beta)$ is an antisymmetric matrix. So $\mbox{Id} + A(\beta)$ is an invertible matrix. Therefore, we obtain:
\begin{equation}\label{1'}
\left\{
\begin{aligned}
&\frac{d\beta}{dt}=(\mbox{Id}+A(\beta))^{-1}B(\beta),\\
&\beta(0)= \left(\int_{\Omega} u_0\cdot\omega^1 \,dx, \int_{\Omega} u_0\cdot\omega^2 \,dx, \cdots, \int_{\Omega} u_0\cdot\omega^N \,dx \right)^T.
\end{aligned}
\right.
\end{equation}
On the other hand, by Lemma 2.1 it is also easy to see that the right side of (\ref{1'}) is locally Lipschitz continuous. By Picard theorem, we known that there is a $\tau>0$ such that the solution of (\ref{1}) exists in $[0,\tau]$. Hence, it follows from (\ref{1}) that
\begin{equation}\label{3}
\left\{
\begin{array}{llll}
  \aligned
  \int_{\Omega}u^{N}_t\cdot\omega^i\,dx = &\,\alpha\int_{\Omega}[\mathfrak{J}(u^N),u^N_t]\omega^i dx+\varepsilon\int_{\Omega}\Delta u^N\omega^i\,dx\\
   &-\int_{\Omega}\left[\mathfrak{J}(u^N),h_d(u^N)\right]\omega^i\,dx-\int_{\Omega}[\mathfrak{J}(u^N),\Delta u^N]\omega^i\,dx\\
   &+\int_{\Omega}[\mathfrak{J}(u^N),(\nabla_u\Phi)\left(\mathfrak{J}(u^N)\right)]\omega^i\,dx,\\
   u^N(0,\cdot)= u_0^N:= &\sum\limits_{i=1}^N\left(\int_{\Omega}u_0\cdot\omega^i\,d\Omega\right)\omega^i,\s\s\s\s\frac{\partial u^N}{\partial\nu}\Big|_{\partial\Omega}=0.
  \endaligned
\end{array}
\right.
\end{equation}
Multiplying the two sides of (\ref{3}) by $\beta_i^N(t)$ and  summing $i$ from $1$ to $N$, and then integrating by parts, we get:
\[
\frac{1}{2}\frac{d}{dt}\int_{\Omega}|u^N|^2\,dx+\varepsilon\int_{\Omega}|\nabla u^N|^2\,dx=0.
\]
This leads to
\begin{equation}\label{4}
\aligned
&\int_{\Omega}|u^N(t)|^2\,dx+2\varepsilon\int_0^t\int_{\Omega}|\nabla u^N|^2\,dxdt=\int_{\Omega}|u_0^N|^2\,dx\leq \int_{\Omega}|u_0|^2\,dx=\mbox{vol}(\Omega),
\endaligned
\end{equation}
where $\mbox{vol}(\Omega)$ denotes the volume of $\Omega$.

Since
\[
\aligned
&\int_{\Omega}|u^N|^2\,dx&=&\int_{\Omega}(\sum\limits_{i=1}^N\beta_i^N(t)\omega^i)(\sum\limits_{j=1}^N\beta_j^N(t)\omega^j)\,dx\\
&&=&\sum\limits_{i=1}^N\sum\limits_{j=1}^N\beta_i^N(t)\beta_j^N(t)\cdot\delta_{ij}=\sum\limits_{i=1}^N|\beta_i^N(t)|^2,
\endaligned
\]
then, for any $T>0$ and any $i$, $\beta_i^N(t)$ can be extended to $[0,T]$. That is to say, $u^N$ can be extended to $[0,T]$.

Next, we want to obtain some uniform a priori estimates on $u^N$ with respect to $N$. First, multiplying the two sides of (\ref{3}) by $-\lambda_i\beta_i^N$ and summing $i$ from 1 to $N$, then we integrate the obtained identity by parts to get:
\begin{eqnarray}\label{5}
 &&-\int_{\Omega}\left[\mathfrak{J}(u^N),\nabla w^N\right]\cdot\Delta u^N-\alpha\int_{\Omega}\left[\mathfrak{J}(u^N),u^N_t\right]\cdot\Delta u^N\\
&=&\frac{1}{2}\frac{d}{dt}\int_{\Omega}|\nabla u^N|^2+\varepsilon\int_{\Omega}|\Delta u^N|^2+\int_{\Omega}\left[\mathfrak{J}(u^N),\nabla_u\Phi\left(\mathfrak{J}(u^N)\right)\right]\cdot\Delta u^N,\nonumber
\end{eqnarray}
where $w^N$ is given by $w^N:=\nabla N\ast u^N$. Multiplying again the two sides of (\ref{3}) by $\frac{d\beta^N_i}{dt}$, summing $i$ from 1 to $N$ and then integrating by parts, we are led to
\begin{eqnarray}\label{6}
&&\int_{\Omega}\left[\mathfrak{J}(u^N),u^N_t\right]\Delta u^N dx+\int_{\Omega}\left[\mathfrak{J}(u^N),\nabla w^N\right]u_t^N\,dx\\
&=&\int_{\Omega}|u^N_t|^2\,dx+\frac{\varepsilon}{2}\frac{d}{dt}\int_{\Omega}|\nabla u^N|^2\,dx-\int_{\Omega}\left[\mathfrak{J}(u^N),\nabla_u\Phi\left(\mathfrak{J}(u^N)\right)\right]u_t^N\,dx.\nonumber
\end{eqnarray}

Now we multiply the two sides of (\ref{6}) by $\alpha$ and then add the two sides of (\ref{5}) to the two sides of (\ref{6}) respectively to obtain:
\begin{eqnarray}\label{7}
&&\alpha\int_{\Omega}|u^N_t|^2\,dx+\frac{1+\alpha\varepsilon}{2}\frac{d}{dt}\int_{\Omega}|\nabla u^N|^2\,dx+\varepsilon\int_{\Omega}|\Delta u^N|^2\,dx\\
&=&-\int_{\Omega}\left[\mathfrak{J}(u^N),\nabla_u\Phi\left(\mathfrak{J}(u^N)\right)\right]\cdot\Delta u^N\,dx +\alpha\int_{\Omega}\left[\mathfrak{J}(u^N),\nabla_u\Phi\left(\mathfrak{J}(u^N)\right)\right]u_t^N\,dx\nonumber\\
&&-\int_{\Omega}\left[\mathfrak{J}(u^N),\nabla w^N\right]\cdot\Delta u^N\,dx+\alpha\int_{\Omega}\left[\mathfrak{J}(u^N),\nabla w^N\right]u_t^N\,dx\nonumber\\
&\leqslant&M_1\cdot\int_{\Omega}|\Delta u^N|\,dx+M_1\cdot\alpha\cdot\int_{\Omega}|u_t^N|\,dx+\int_{\Omega}|\nabla w^N|\cdot|\Delta u^N|\,dx\nonumber\\
&&+\alpha\cdot\int_{\Omega}|\nabla w^N|\cdot|u_t^N|\,dx\nonumber\\
&\leqslant&\frac{\varepsilon}{2}\int_{\Omega}|\Delta u^N|^2 +\frac{\alpha}{2}\int_{\Omega}|u^N_t|^2\,dx+\frac{(1+\alpha)}{\varepsilon}\cdot\int_{\mathbb{R}^n}|\nabla w^N|^2+(\alpha+\frac{1}{\varepsilon}) M_1^2\mbox{vol}(\Omega),\nonumber
\end{eqnarray}
where $M_1$ depends on the value of $\nabla_u\Phi$ which is restricted in the unit closed ball of the Lie algebra $\mathfrak{g}$. By rearranging the above inequality we can derive
\begin{eqnarray}\label{add}
&&\alpha\int_0^t\int_{\Omega}|u^N_t|^2+(1+\alpha\varepsilon)\cdot\int_{\Omega}|\nabla u^N|^2+\varepsilon\int_0^t\int_{\Omega}|\Delta u^N|^2\nonumber\\
&\leqslant&(2\alpha M_1^2+\frac{2M^2_1}{\varepsilon})\mbox{vol}(\Omega)t+\frac{2(1+\alpha)}{\varepsilon}\int_0^t\int_{\mathbb{R}^n}|\nabla w^N|^2+(1+\alpha\varepsilon)\int_{\Omega}|\nabla u_0|^2,
\end{eqnarray}
where we have used the fact that $$\int_{\Omega}|\nabla u^N_0|^2\leqslant\int_{\Omega}|\nabla u_0|^2.$$

Lemma \ref{lem} tells us that there hold
\begin{eqnarray}\label{ZL3}
\int_{\mathbb{R}^n}|\nabla w^N|^2\,dx\leqslant\int_{\Omega}|u^N|^2\,dx
\end{eqnarray}
for any $N\geqslant0$, and
\begin{eqnarray}\label{ZL4}
\int_{\mathbb{R}^n}|\nabla w^{N_1}-\nabla w^{N_2}|^2\,dx\leqslant\int_{\Omega}|u^{N_1}-u^{N_2}|^2\,dx
\end{eqnarray}
for any $N_1, N_2\geqslant0$.

In view of $(\ref{4})$ and $(\ref{ZL3})$ we can see from (\ref{add}) that
\begin{eqnarray}\label{ZL6}
&&\alpha\int_0^t\int_{\Omega}|u^N_t|^2+(1+\alpha\varepsilon)\cdot\int_{\Omega}|\nabla u^N|^2+\varepsilon\int_0^t\int_{\Omega}|\Delta u^N|^2\\
&\leqslant&\left(2\alpha M_1^2+\frac{2M^2_1+2(1+\alpha)}{\varepsilon}\right)\mbox{vol}(\Omega)t+(1+\alpha\varepsilon)\int_{\Omega}|\nabla u_0|^2.\nonumber
\end{eqnarray}
Hence it is easy to conclude
\begin{lem}\label{lem1}
The approximate solution sequence $\{u^N\}$ to (\ref{1}) satisfies

$\bullet$ $\{u^N\}$ is a bounded sequence in $L^{\infty}([0,T],H^1(\Omega,\mathfrak{g}))$;

\medskip
$\bullet$ $\{u^N_t\}$ is a bounded sequence in $L^{2}([0,T],L^2(\Omega,\mathfrak{g}))$;

\medskip
$\bullet$ $\{\Delta u^N\}$ is a bounded sequence in $L^{2}([0,T],L^2(\Omega,\mathfrak{g}))$;

\medskip
$\bullet$ $\{\nabla u^N\}$ is a bounded sequence in $L^{\infty}([0,T],L^2(\Omega,\mathbb{R}^n\otimes\mathfrak{g}))$.
\end{lem}

\medskip
By the property of weak limits and Aubin-Lions Lemma, from Lemma \ref{lem1} we deduce that there exists a $v^{\varepsilon}\in W^{2,1}_2(\Omega\times[0,T],\mathfrak{g})$ and a subsequence of $\{u^N\}$ which is also denoted by $\{u^N\}$ such that

$\bullet$ $u^N\rightharpoonup v^{\varepsilon}\,\,\,\,\mbox{weakly* in}\,\,L^{\infty}([0,T],H^1(\Omega,\mathfrak{g}))$;

\medskip
$\bullet$ $u^N\rightarrow v^{\varepsilon}\,\,\,\,\mbox{strongly in}\,\,L^{\infty}([0,T],L^2(\Omega,\mathfrak{g}))$;

\medskip
$\bullet$ $u^N\rightarrow v^{\varepsilon}\,\,\,\,a.e.\,\,\Omega\times[0,T]$;

\medskip
$\bullet$ $u^N_t\rightharpoonup v^{\varepsilon}_t\,\,\,\,\mbox{weakly in}\,\,L^{2}([0,T],L^2(\Omega,\mathfrak{g}))$;

\medskip
$\bullet$ $\Delta u^N\rightharpoonup \Delta v^{\varepsilon}\,\,\,\,\mbox{weakly in}\,\,L^2([0,T],L^2(\Omega,\mathfrak{g}))$;

\medskip
$\bullet$ $\nabla u^N\rightharpoonup \nabla v^{\varepsilon}\,\,\,\,\mbox{weakly* in}\,\,L^{\infty}([0,T],L^2(\Omega,\mathbb{R}^n\otimes\mathfrak{g}))$;

\medskip
$\bullet$ $h_d(u^N)=-\nabla w^N\rightarrow h_d(v^{\varepsilon})\,\,\,\,\mbox{strongly in}\,\,L^{\infty}([0,T],L^2(\mathbb{R}^n,\mathfrak{g}))$, where we have used inequality $(\ref{ZL4})$;

\medskip
$\bullet$ $h_d(u^N)=-\nabla w^N\rightarrow h_d(v^{\varepsilon})\,\,\,\,a.e.\,\,\mathbb{R}^n\times[0,T]$.
\medskip

From the following facts
$$
||u^N||_{L^{\infty}([0,T],H^1(\Omega,\mathfrak{g}))}\leq C_{12}(\varepsilon)
$$
and $u^N\rightharpoonup v^{\varepsilon}\,\,\,\,weakly*\,\,\,\,in\,\,L^{\infty}([0,T],H^1(\Omega,\mathfrak{g}))$, we have
\begin{eqnarray}\label{ZL5}
||v^{\varepsilon}||_{L^{\infty}([0,T],H^1(\Omega,\mathfrak{g}))}\leq C_{12}(\varepsilon).
\end{eqnarray}
By the same method as we prove $(\ref{ZL5})$, from (\ref{ZL6}) it follows that
\[
\alpha\int_0^T\int_{\Omega}|v_t^{\varepsilon}|^2\,dx\,dt\leq C_{13}(\varepsilon).
\]
Hence, it is easy to see that
\[
\left[\mathfrak{J}(u^N),u^N_t\right]\rightharpoonup\left[\mathfrak{J}(v^{\varepsilon}),v^{\varepsilon}_t\right]
\]
weakly in $L^{2}([0,T],L^2(\Omega,\mathfrak{g}))$, and
\[
\left[\mathfrak{J}(u^N),\Delta u^N\right]\rightharpoonup\left[\mathfrak{J}(v^{\varepsilon}),\Delta v^{\varepsilon}\right]
\]
weakly in $L^{2}([0,T],L^2(\Omega,\mathfrak{g}))$.

Fixing $r\in \mathbb{Z}^+$ and taking any $N\geq r$, we multiply two sides of (\ref{3}) by $\eta^i(t)$ which belongs to $C^{\infty}([0,T],\mathfrak{g})$ and sum $i$ from 1 to $r$, and then integrate the obtained identity on $[0,T]$ to derive
\[
  \aligned
&\int_0^Tdt\int_{\Omega}u^{N}_t\cdot\Phi^r\,dx=\alpha\int_0^Tdt\int_{\Omega}[\mathfrak{J}(u^N),u^N_t]\Phi^r dx+\varepsilon\int_0^Tdt\int_{\Omega}\Delta u^N\Phi^r\,dx\\
&-\int_0^Tdt\int_{\Omega}\left[\mathfrak{J}(u^N),h_d(u^N)\right]\Phi^r\,dx-\int_0^Tdt\int_{\Omega}[\mathfrak{J}(u^N),\Delta u^N]\Phi^r\,dx\\
&+\int_0^Tdt\int_{\Omega}[\mathfrak{J}(u^N),(\nabla_u\Phi)\left(\mathfrak{J}(u^N)\right)]\Phi^r\,dx,
  \endaligned
\]
where
\[
\Phi^r(x,t)=\sum\limits_{i=1}^r\omega^i(x)\eta^i(t).
\]
Letting $N$ tends to $\infty$ in the above identity, we get:
\[
  \aligned
\int_0^Tdt\int_{\Omega}v^{\varepsilon}_t\cdot\Phi^r\,dx=&\alpha\int_0^Tdt\int_{\Omega}[\mathfrak{J}(v^{\varepsilon}),v^{\varepsilon}_t]\Phi^r dx+\varepsilon\int_0^Tdt\int_{\Omega}\Delta v^{\varepsilon}\Phi^r\,dx\\
&-\int_0^Tdt\int_{\Omega}\left[\mathfrak{J}(v^{\varepsilon}),h_d(v^{\varepsilon})\right]\Phi^r\,dx-\int_0^Tdt\int_{\Omega}[\mathfrak{J}(v^{\varepsilon}),\Delta v^{\varepsilon}]\Phi^r\,dx\\
&+\int_0^Tdt\int_{\Omega}[\mathfrak{J}(v^{\varepsilon}),(\nabla_u\Phi)\left(\mathfrak{J}(v^{\varepsilon})\right)]\Phi^r\,dx.
  \endaligned
\]
Since the functions, which are of type $\Phi^r(x,t)$, are dense in $L^2([0,T],L^2(\Omega,\mathfrak{g}))$, we known that, in the sense of distribution, there holds
\begin{equation}\label {Key'}
\left\{ \begin{aligned}
& v^{\varepsilon}_t = \varepsilon\Delta v^{\varepsilon}+[\mathfrak{J}(v^{\varepsilon}),\,\alpha\partial_t v^{\varepsilon}-\Delta v^{\varepsilon}-h_d(v^{\varepsilon})+\nabla_u\Phi(\mathfrak{J}(v^{\varepsilon}))],\\
& v^{\varepsilon}(\cdot,0)=u_0:\Omega\longrightarrow S_{\mathfrak{g}}(1).
\end{aligned} \right.
\end{equation}

 Next, we would like to show that $\frac{\partial v^{\varepsilon}}{\partial\nu}\big|_{\partial\Omega}=0$. Indeed, for any $\phi\in C^{\infty}(\Omega)$, we have
$$\int_{\Omega}\Delta u^N\cdot\phi+\sum\limits^n_{p=1}\int_{\Omega}\nabla_p u^N\cdot\nabla_p\phi=0,$$
since $\frac{\partial u^N}{\partial\nu}\big|_{\partial\Omega}=0$. Letting $N$ tends to $\infty$ yields
$$\int_{\Omega}\Delta v^{\varepsilon}\cdot\phi+\sum\limits^n_{p=1}\int_{\Omega}\nabla_p v^{\varepsilon}\cdot\nabla_p\phi=0.$$
The arbitrariness of $\phi$ implies $\frac{\partial v^{\varepsilon}}{\partial\nu}\Big|_{\partial\Omega}=0$.

\subsection{Derivation of the limit equation and some uniform estimates}

\medskip
Choosing $$v^{\varepsilon}-v^{\varepsilon}\frac{min\{1,|v^{\varepsilon}|\}}{|v^{\varepsilon}|}$$ as a test function of the above equation (\ref{Key'}), we obtain:
\[
\int_{\Omega}v^{\varepsilon}_t\cdot\left(v^{\varepsilon}-v^{\varepsilon}\frac{min\{1,|v^{\varepsilon}|\}}{|v^{\varepsilon}|}\right)\,dx=\varepsilon\int_{\Omega}\Delta v^{\varepsilon}\cdot\left(v^{\varepsilon}-v^{\varepsilon}\frac{min\{1,|v^{\varepsilon}|\}}{|v^{\varepsilon}|}\right)\,dx.
\]
By a simple computation we can see from the above identity
\begin{equation}\label{Key1'}
\aligned
&\frac{1}{2}\frac{d}{dt}\int_{|v^{\varepsilon}|> 1}|v^{\varepsilon}|^2\left(1-\frac{1}{|v^{\varepsilon}|}\right)\,dx+\varepsilon\int_{|v^{\varepsilon}|> 1}\frac{|\nabla v^{\varepsilon}\cdot v^{\varepsilon}|^2}{|v^{\varepsilon}|^3}\,dx\\
=&\frac{1}{2}\int_{|v^{\varepsilon}|>1}\frac{v^{\varepsilon}\cdot v^{\varepsilon}_t}{|v^{\varepsilon}|}\,dx-\varepsilon\int_{|v^{\varepsilon}|>1}|\nabla v^{\varepsilon}|^2\left(1-\frac{1}{|v^{\varepsilon}|}\right)\,dx.
\endaligned
\end{equation}
Taking $$\frac{v^{\varepsilon}(max\{|v^{\varepsilon}|,1\}-1)}{|v^{\varepsilon}|(|v^{\varepsilon}|-1+\delta)}$$ as another test function of (\ref{Key'}), we get:
\[
\aligned
&\int_{|v^{\varepsilon}|>1}\frac{v^{\varepsilon}\cdot v^{\varepsilon}_t}{|v^{\varepsilon}|}\cdot\frac{|v^{\varepsilon}|-1}{|v^{\varepsilon}|-1+\delta}\,dx\\
=&-\varepsilon\int_{\Omega}\nabla v^{\varepsilon}\cdot\nabla\Big[ \frac{v^{\varepsilon}(\max\{|v^{\varepsilon}|,1\}-1)}{|v^{\varepsilon}|(|v^{\varepsilon}|-1+\delta)} \Big]\,dx\\
=&-\varepsilon\int_{|v^{\varepsilon}|>1}|\nabla v^{\varepsilon}|^2\frac{|v^{\varepsilon}|-1}{|v^{\varepsilon}|(|v^{\varepsilon}|-1+\delta)}\,dx\\
&-\varepsilon\int_{|v^{\varepsilon}|>1}\frac{|\nabla v^{\varepsilon}\cdot v^{\varepsilon}|^2}{|v^{\varepsilon}|}\cdot\frac{-|v^{\varepsilon}|^2+2|v^{\varepsilon}|-1+\delta}{|v^{\varepsilon}|^2(|v^{\varepsilon}|+\delta-1)^2}\,dx.
\endaligned
\]
By the dominated convergence theorem, letting $\delta\rightarrow0$ we derive from the above
\begin{equation}\label{Key2'}
\int_{|v^{\varepsilon}|>1}\frac{v^{\varepsilon}\cdot v^{\varepsilon}_t}{|v^{\varepsilon}|}\,dx=-\varepsilon\int_{|v^{\varepsilon}|>1}\frac{|\nabla v^{\varepsilon}|^2}{|v^{\varepsilon}|}\,dx+\varepsilon\int_{|v^{\varepsilon}|>1}\frac{|\nabla v^{\varepsilon}\cdot v^{\varepsilon}|^2}{|v^{\varepsilon}|^3}\,dx.
\end{equation}
Combining (\ref{Key1'}) and (\ref{Key2'}) yields
\[
\frac{d}{dt}\int_{|v^{\varepsilon}|>1}|v^{\varepsilon}|^2\Big(1-\frac{1}{|v^{\varepsilon}|}\Big)\,dx\leq 0.
\]
This means that the following function
\[
q(t):=\int_{|v^{\varepsilon}|>1}|v^{\varepsilon}|^2\Big(1-\frac{1}{|v^{\varepsilon}|}\Big)\,dx
\]
is decreasing non-negative function. Noting $|v^{\varepsilon}(\cdot,0)|=|u_0|=1$, i.e. $q(0)=0$, we can see that $q(t)\equiv0$ for any $t>0$. Therefore, we have $|v^{\varepsilon}|\leq1$. Hence, we have shown the following lemma:

\begin{lem}
For any fixing $\varepsilon>0$, the auxiliary approximation equation (\ref{A0}) admits a weak solution $v^{\varepsilon}$ belonging to $W^{2,1}_2(\Omega\times[0, T],\mathfrak{g})$, which satisfies that for any $t\in [0, T]$
$$|v^{\varepsilon}(t)|\leq 1,\s\s a.e. \s \Omega.$$
\end{lem}

Immediately it follows from the above lemma that equation (\ref{Key'}) becomes into the following
\begin{equation}\label{Key3'}
\partial_t v^{\varepsilon}  = \varepsilon\Delta v^{\varepsilon}+[v^{\varepsilon},\,\alpha\partial_t v^{\varepsilon}-\Delta v^{\varepsilon}-h_d(v^{\varepsilon})+\nabla_u\Phi(v^{\varepsilon})],
\end{equation}
with initial value $v^{\varepsilon}(\cdot, 0)=u_0$ and $\frac{\partial v^{\varepsilon}}{\partial\nu}\Big|_{\partial\Omega}=0$.

Now, multiplying the two sides of $(\ref{Key3'})$ by $v^{\varepsilon}$ and integrating it on $\Omega\times[0,t]$, we get:
\begin{equation}\label{8}
\int_{\Omega}(|v^{\varepsilon}(t)|^2-1)\,dx+2\varepsilon\int_0^t\int_{\Omega}|\nabla v^{\varepsilon}|^2\,dxd\tau=0.
\end{equation}
Multiplying the both sides of $(\ref{Key3'})$ by $\Delta v^{\varepsilon}$ and integrating it on $\Omega\times[0,t]$ lead to
\begin{equation*}
\int_0^t\int_{\Omega} v^{\varepsilon}_t\cdot\Delta v^{\varepsilon}  =\varepsilon\int_0^t\int_{\Omega}|\Delta v^{\varepsilon}|^2+\int_0^t\int_{\Omega}[v^{\varepsilon},\,\alpha\partial_t v^{\varepsilon}-h_d(v^{\varepsilon})+\nabla_u\Phi(v^{\varepsilon})]\cdot\Delta v^{\varepsilon},
\end{equation*}
which means
\begin{eqnarray}\label{9}
&&\frac{1}{2}\int_{\Omega}|\nabla u_0|^2+\int_0^t\int_{\Omega}[v^{\varepsilon},\nabla^2w^{\varepsilon}]\cdot\nabla v^{\varepsilon}+\int_0^t\int_{\Omega}[v^{\varepsilon},\nabla^2_u\Phi(v^{\varepsilon})\nabla v^{\varepsilon}]\cdot\nabla v^{\varepsilon}\\
&=&\frac{1}{2}\int_{\Omega}|\nabla v^{\varepsilon}(t)|^2+\varepsilon\int_0^t\int_{\Omega}|\Delta v^{\varepsilon}|^2+\alpha\int_0^t\int_{\Omega}[v^{\varepsilon},\,\partial_t v^{\varepsilon}]\Delta v^{\varepsilon},\nonumber
\end{eqnarray}
where $w^{\varepsilon}$ is defined as $w^{\varepsilon}:=\nabla N\ast v^{\varepsilon}$.
Multiplying again the two sides of $(\ref{Key3'})$ by $v_t^{\varepsilon}$, integrating it on $\Omega\times[0,t]$ and then integrating by parts we have
\begin{eqnarray}\label{10}
&&\int_0^t\int_{\Omega}[v^{\varepsilon},\nabla w^{\varepsilon}]\cdot v_t^{\varepsilon}+\frac{\varepsilon}{2}\int_{\Omega}|\nabla u_0|^2+\int_0^t\int_{\Omega}[v^{\varepsilon},\nabla_u\Phi(v^{\varepsilon})]\cdot v_t^{\varepsilon}\\
&=&\frac{\varepsilon}{2}\int_{\Omega}|\nabla v^{\varepsilon}(t)|^2+\int_0^t\int_{\Omega}|v_t^{\varepsilon}|^2-\int_0^t\int_{\Omega}[v^{\varepsilon},v_t^{\varepsilon}]\cdot\Delta v^{\varepsilon}.\nonumber
\end{eqnarray}

Then, by multiplying the two sides of (\ref{10}) by $\alpha$ and then adding respectively the two sides of (\ref{9}) to the two sides of (\ref{10}), we obtain:
\begin{eqnarray}\label{11}
&&\frac{1+\alpha\varepsilon}{2}\int_{\Omega}|\nabla v^{\varepsilon}(t)|^2+\alpha\int_0^t\int_{\Omega}|v_t^{\varepsilon}|^2+\varepsilon\int_0^t\int_{\Omega}|\Delta v^{\varepsilon}|^2\\
&=&\alpha\int_0^t\int_{\Omega}[v^{\varepsilon},\nabla w^{\varepsilon}]\cdot v_t^{\varepsilon}+\frac{1+\alpha\varepsilon}{2}\int_{\Omega}|\nabla u_0|^2+\alpha\int_0^t\int_{\Omega}[v^{\varepsilon},\nabla_u\Phi(v^{\varepsilon})]\cdot v_t^{\varepsilon}\nonumber\\
&&+\int_0^t\int_{\Omega}[v^{\varepsilon},\nabla^2w^{\varepsilon}]\cdot\nabla v^{\varepsilon}+\int_0^t\int_{\Omega}[v^{\varepsilon},\nabla^2_u\Phi(v^{\varepsilon})\nabla v^{\varepsilon}]\cdot\nabla v^{\varepsilon}\nonumber\\
&\leqslant&\alpha \int_0^t\int_{\Omega}|\nabla w^{\varepsilon}|\cdot|v_t^{\varepsilon}|+\frac{1+\alpha\varepsilon}{2}\int_{\Omega}|\nabla u_0|^2+\alpha M_1\int_0^t\int_{\Omega}|v_t^{\varepsilon}|\nonumber\\
&&+\int_0^t\int_{\Omega}|\nabla^2w^{\varepsilon}|\cdot|\nabla v^{\varepsilon}|+M_2\int_0^t\int_{\Omega}|\nabla v^{\varepsilon}|^2\nonumber\\
&\leqslant&\frac{\alpha}{2}\int_0^t\int_{\Omega}|v_t^{\varepsilon}|^2+\alpha\int_0^t\int_{\mathbb{R}^n}|\nabla w^{\varepsilon}|^2+\frac{1+\alpha\varepsilon}{2}\int_{\Omega}|\nabla u_0|^2\nonumber\\
&&+(M_2+1/2)\int_0^t\int_{\Omega}|\nabla v^{\varepsilon}|^2+\frac{1}{2}\int_0^t\int_{\Omega}|\nabla^2w^{\varepsilon}|^2+\alpha\cdot M_1^2\cdot\mbox{vol}(\Omega)\cdot t\nonumber\\
&\leqslant&\frac{\alpha}{2}\int_0^t\int_{\Omega}|v_t^{\varepsilon}|^2+\frac{1+\alpha\varepsilon}{2}\int_{\Omega}|\nabla u_0|^2+\alpha\cdot (1+M_1^2)\cdot\mbox{vol}(\Omega)\cdot t\nonumber\\
&&+\left(M_2+\frac{1+C}{2}\right)\int_0^t\int_{\Omega}|\nabla v^{\varepsilon}|^2+\frac{C}{2}\mbox{vol}(\Omega)t,\nonumber
\end{eqnarray}
where $M_2$ depends on the value of $\nabla^2_u\Phi$, which is restricted on the unit closed ball of the Lie algebra $\mathfrak{g}$, and we have used the following facts (see Lemma 2.1)
$$\int_{\mathbb{R}^n}|\nabla w^{\varepsilon}|^2\leqslant\int_{\Omega}|v^{\varepsilon}|^2\leqslant\int_{\Omega}1=\mbox{vol}(\Omega)$$
and
\begin{eqnarray}\label{ZL7}
\int_{\Omega}|\nabla^2 w^{\varepsilon}|^2\leqslant C\|v^{\varepsilon}\|^2_{W^{1,2}(\Omega)}.
\end{eqnarray}

It follows (\ref{11}) that
\begin{eqnarray}\label{12}
&&\frac{1+\alpha\varepsilon}{2}\int_{\Omega}|\nabla v^{\varepsilon}(t)|^2+\frac{\alpha}{2}\int_0^t\int_{\Omega}|v_t^{\varepsilon}|^2+\varepsilon\int_0^t\int_{\Omega}|\Delta v^{\varepsilon}|^2\\
&\leqslant&\frac{1+\alpha\varepsilon}{2}\int_{\Omega}|\nabla u_0|^2+(M_2+\frac{1+C}{2})\int_0^t\int_{\Omega}|\nabla v^{\varepsilon}|^2 +(\alpha (1+M_1^2)+\frac{C}{2})\mbox{vol}(\Omega)\cdot t.\nonumber
\end{eqnarray}
Hence, Gronwall Inequality tells us that, for any $t\in [0, T]$, there holds
$$\int_{\Omega}|\nabla v^{\varepsilon}(t)|^2\leqslant M_3,$$
where $M_3$ depends upon $T, M_1, M_2, \Omega, \alpha$ but not upon $\varepsilon$. Substituting the upper bound of the above quantity into $(\ref{12})$ yields
\begin{eqnarray}\label{13}
&&(1+\alpha\varepsilon)\int_{\Omega}|\nabla v^{\varepsilon}(t)|^2+\alpha\int_0^t\int_{\Omega}|v_t^{\varepsilon}|^2+2\varepsilon\int_0^t\int_{\Omega}|\Delta v^{\varepsilon}|^2\\
&\leqslant&(1+\alpha\varepsilon)\int_{\Omega}|\nabla u_0|^2+2\alpha\cdot (1+M_1^2)\cdot\mbox{vol}(\Omega)\cdot t+(2M_2+1+n)M_3\cdot t.\nonumber
\end{eqnarray}
This implies that, for any $t\leq T$, there holds true
$$\alpha\int_0^t\int_{\Omega}|v_t^{\varepsilon}|^2\leq M_4.$$
Thus we have established the following
\begin{lem}\label{lem3}
(1). In the case $\alpha=0$, for any $T>0$ there holds true for the solution $v^{\varepsilon}$ to (\ref{Key3'}) belonging to $L^\infty_{loc}(\mathbb{R}^+, H^1(\mathbb{T},\mathfrak{g}))$
$$\sup_{t\in[0, T]}\int_{\Omega}|\nabla v^{\varepsilon}(t)|^2\leqslant M_3(T),$$
where $M_3$ does not depend on $\varepsilon$.

\noindent(2). In the case $\alpha>0$, for any $T>0$ there holds true for the solution $v^{\varepsilon}$ to (\ref{Key3'}) belonging to $W^{2,1}_2(\Omega\times [0, T],\mathfrak{g})$
$$\sup_{t\in[0, T]}\int_{\Omega}|\nabla v^{\varepsilon}(t)|^2\leqslant M_3(T),\s\s\s\s \alpha\int_0^T\int_{\Omega}|v_t^{\varepsilon}|^2\leq M_4(T),$$
where $M_3$ and $M_4$ do not depend on $\varepsilon$.
\end{lem}

\subsection{Proof of Theorem \ref{thm}}

Now we return to present the proof of of Theorem \ref{thm} and need to consider the following two cases:

\textbf{Case 1}: $\alpha>0$. From Lemma \ref{lem3} we know that $\{v^{\varepsilon}\}$ is a bounded sequence in the Sobolev space $L^{\infty}([0,T],H^1(\Omega,\mathfrak{g}))$ and $\{v_t^{\varepsilon}\}$ is a bounded sequence in $L^2([0,T],L^2(\Omega,\mathfrak{g}))$. So, by the property of weak limits and Aubin-Lions Lemma, there is a $u$ and a subsequence of $\{v^{\varepsilon}\}$ which is also denoted by $\{v^{\varepsilon}\}$ such that:

$\bullet$ $v^{\varepsilon}\rightharpoonup u$ weakly $*$ in $L^{\infty}([0,T],H^1(\Omega,\mathfrak{g}))$;
\medskip

$\bullet$ $v^{\varepsilon}\rightarrow u$ strongly in $L^{\infty}([0,T],L^2(\Omega,\mathfrak{g}))$;
\medskip

$\bullet$ $h_d(v^{\varepsilon})\rightarrow h_d(u)$ strongly in $L^{\infty}([0,T],L^2(\Omega,\mathfrak{g}))$;
\medskip

$\bullet$ $v^{\varepsilon}\rightarrow u$ a.e. $\Omega\times[0,T]$;
\medskip

$\bullet$ $v^{\varepsilon}_t\rightharpoonup u_t$ weakly in $L^{2}([0,T],L^2(\Omega,\mathfrak{g}))$.\\
\noindent Here we have used the fact that
$$\int_{\mathbb{R}^n}|h_d(v^{\varepsilon_i})-h_d(v^{\varepsilon_j})|^2\leqslant\int_{\Omega}|v^{\varepsilon_i}-v^{\varepsilon_j}|^2.$$
Therefore, letting $\varepsilon$ in (\ref{8}) tends to 0, we have:
\[
\int_{\mathbb{T}}(|u|^2-1)\,d\mathbb{T}=0.
\]
This leads to $|u|=1$ for $a.e.\, x\in\mathbb{T}$ for all $t\in[0,T]$ since $|v^{\varepsilon}|\leq1$ implies that $|u|\leq1$.

For any $\varphi\in C^{\infty}([0,T]\times\mathbb{T},\mathfrak{g})$, we have
\begin{eqnarray*}
&&\int_0^Tdt\int_{\Omega}\partial_tv^{\varepsilon}\cdot\varphi\,d\Omega
-\alpha\cdot\int_0^Tdt\int_{\Omega}[v^{\varepsilon},\partial_tv^{\varepsilon}]\cdot\varphi\,d\Omega+\varepsilon\sum\limits_{p=1}^n\int_0^Tdt\int_{\Omega}\nabla_p v^{\varepsilon}\cdot\nabla_p\varphi\\
&=&\sum\limits_{p=1}^n\int_0^Tdt\int_{\Omega}[v^{\varepsilon},\nabla_{p}v^{\varepsilon}]\cdot\nabla_{p}\varphi\,d\Omega-\int_0^Tdt\int_{\Omega}[v^{\varepsilon},h_d(v^{\varepsilon})-\nabla_u\Phi(v^{\varepsilon})]\cdot\varphi\,d\Omega.
\end{eqnarray*}
By letting $\varepsilon$ tend to $0$ in the above identity, it is easy to see
\begin{eqnarray*}
&&\int_0^Tdt\int_{\Omega}\partial_tu\cdot\varphi\,d\Omega-\alpha\cdot\int_0^Tdt\int_{\Omega}[u,\,\partial_tu]\cdot\varphi\,d\Omega \\
&=&\sum\limits_{p=1}^n\int_0^T dt\int_{\Omega}[u,\,\nabla_{p}u]\cdot\nabla_{p}\varphi\,d\Omega-\int_0^Tdt\int_{\Omega}[u,\,h_d(u)-\nabla_u\Phi(u)]\cdot\varphi\,d\Omega.
\end{eqnarray*}
This means that $u$ is a weak solution to (\ref{spin:1}) with $\alpha>0$.

\medskip
\textbf{Case 2:} $\alpha=0$. We have known that, as $\alpha>0$, $\{v^{\varepsilon}\}$ is a bounded sequence in the Sobolev space $L^{\infty}([0,T],H^1(\Omega,\mathfrak{g}))\cap W^{1,1}_2(\Omega\times[0,T], \mathfrak{g})$. By the same argument as in the above, letting $\varepsilon$ in (\ref{8}) tends to $0$ and denoting the limit of $v^\varepsilon$ by $u^\alpha$ we conclude that, for any $\varphi\in C^{\infty}([0,T]\times\mathbb{T},\mathfrak{g})$,
\begin{eqnarray}\label{fin}
&&\int_{\Omega}(u^\alpha(T)\cdot\varphi(T)-u_0\cdot\varphi(0))-\int_0^Tdt\int_{\Omega}u^\alpha\cdot\partial_t\varphi
-\alpha\int_0^Tdt\int_{\Omega}[u^\alpha,\,\partial_tu^\alpha]\cdot\varphi\\ \nonumber
&=&\sum\limits_{p=1}^n\int_0^Tdt\int_{\Omega}[u^\alpha,\nabla_{p}u^\alpha]\cdot\nabla_{p}\varphi\,d\Omega
-\int_0^Tdt\int_{\Omega}[u^\alpha,\, h_d(u^\alpha)-\nabla_u\Phi(u^\alpha)]\cdot\varphi\,d\Omega.
\end{eqnarray}

Obviously, we have $u^\alpha$ is uniformly bounded in $L^{\infty}([0,T],H^1(\Omega,\mathfrak{g}))$. Therefore, there exists a $u\in L^{\infty}([0,T],H^1(\Omega,\mathfrak{g}))$ and a subsequence of $u^{\alpha}$, which is still denoted by $u^{\alpha}$, such that $u^\alpha\to u$ weakly $\ast$ in $L^{\infty}([0,T],H^1(\Omega,\mathfrak{g}))$. Hence, it is easy to see that $u^\alpha\to u$ a.e. $\Omega\times[0,T]$. Moreover, as $h_d(u^\alpha)=-\nabla w^\alpha$, by the regularity theory of elliptic equation and Lemma 2.1 we also have that $h_d(u^\alpha)\to h_d(u)$ strongly in the space $L^2(\Omega, \mathfrak{g})$ and $h_d(u^\alpha)\to h_d(u)$ a.e. on $\mathbb{R}^n\times[0,T]$.

Noting that, as $\alpha\rightarrow 0$,
$$\alpha\int_0^T\int_{\Omega}|\partial_t u^\alpha|^2\leq M_4$$
implies
$$\alpha\cdot\int_0^Tdt\int_{\Omega}[u^\alpha,\,\partial_tu^\alpha]\cdot\varphi\,d\Omega\longrightarrow 0,$$
we let $\alpha\to 0$ in the above (\ref{fin}) to derive
\begin{eqnarray*}
&&\int_{\Omega}(u(T)\cdot\varphi(T)-u_0\cdot\varphi(0))\,d\Omega + \int_0^Tdt\int_{\Omega}[u,\,h_d(u)-\nabla_u\Phi(u)]\cdot\varphi\,d\Omega\\
&=&\int_0^Tdt\int_{\Omega}u\cdot\partial_t\varphi\,d\Omega+ \sum\limits_{p=1}^n\int_0^T dt\int_{\Omega}[u,\,\nabla_{p}u]\cdot\nabla_{p}\varphi\,d\Omega.
\end{eqnarray*}
By the definition we know that $u$ is a weak solution to (\ref{spin:1}) with $\alpha=0$. Thus, the proof of Theorem \ref{thm} finishes.
\endproof

\begin{rem}
It is worth to point out that, in the case $\alpha=0$, one can also prove the theorem directly. Indeed, from the above arguments we can see easily that the proof of the case $\alpha=0$ goes almost the same as the proof of Theorem \ref{thm} except for one needs to make $\alpha=0$ in the auxiliary approximation equation.
\end{rem}

\begin{rem}
We can also consider the following Landau-Lifshitz-Gilbert system with Lie algebra value:
\begin{equation*}
\left\{ \begin{aligned}
&\partial_tu-\alpha[u,\,\partial_tu]=-[u,\,\Delta u+h_d-\nabla_u\Phi] + F(x,t,u),              & x\in\Omega,\\
& u(\cdot,0)=u_0:\Omega\longrightarrow S_{\mathfrak{g}}(1), \s\s \frac{\partial u}{\partial \nu}\big|_{\partial\Omega}=0.
\end{aligned} \right.
\end{equation*}
Here, $\langle F(x, t, z),\, z\rangle\equiv 0$. In fact, from the above arguments in the section and the following section we can see that it is not difficult to address the existence of the weak solution to the above LLG system.
\end{rem}

\section{Generalized Inhomogeneous LLG Equations on Closed Manifolds}\label{section2}
In this section we will show the well-posedness of global solutions to (\ref{lie:1}). If $f(x)$ is smooth enough and $\dim(\mathbb{T})$ is not $2$, $f(x)$ can be absorbed by a conformal transformation of the metric $h$ on $\mathbb{T}$. Here, we focus on the case $\dim(\mathbb{T})=2$. Recall that the equation we consider is as follows
\begin{eqnarray}\label{Popu}
\left\{
\begin{array}{ll}
  \aligned
  \alpha_0v_t+\alpha\left[v,\, v_t\right]=&\,\left[v,\, f\Delta v+\nabla f\cdot\nabla v\right]+F\left(x,t,v\right),\\
  v(0,\cdot)=u_0.\s\s\s\,\, &
  \endaligned
\end{array}
\right.
\end{eqnarray}
Here, $\alpha_0>0$ and $\alpha\geq 0$ are two constants.

As before, We still need to employ an auxiliary approximation equation as follows
\[
\left\{
\begin{array}{llll}
  \aligned
  &\alpha_0v^{\varepsilon}_t+\alpha\bigg[\frac{v^{\varepsilon}}{max\{|v^{\varepsilon}|,1\}},v^{\varepsilon}_t\bigg]=\varepsilon(f\Delta v^{\varepsilon}+\nabla f\cdot\nabla v^{\varepsilon})\\
  &\quad\quad+\left[\frac{v^{\varepsilon}}{max\{|v^{\varepsilon}|,1\}},f\Delta v^{\varepsilon}+\nabla f\cdot\nabla v^{\varepsilon}\right]+F\left(x,t,\frac{v^{\varepsilon}}{max\{|v^{\varepsilon}|,1\}}\right),\\
  &v^{\varepsilon}(0,\cdot)=u_0.
  \endaligned
\end{array}
\right.
\]
Here, $F\left(x,t,z\right)$ has been extended on $\mathbb{T}\times\mathbb{R}^+\times \overline{B} _{\mathfrak{g}}(1)$ where $\overline{B} _{\mathfrak{g}}(1)$ is a unit closed ball in $\mathfrak{g}$. In fact, we set
\begin{equation*}
\tilde{F}(x, t, z)=
\left\{
\begin{aligned}
&\zeta(|z|^2)F\left(x,,t, \frac{z}{\max\{\delta_0,|z|\}}\right),\s\s& |z|^2>\delta_0,\\
&0, \s\s &|z|^2\leq\delta_0,
\end{aligned}\right.
\end{equation*}
where $\zeta(t):[0, 1]\to [0, 1]$ is a smooth function with $\zeta(t)\equiv 0$ on $[0, 2\delta_0]$ ($2\delta_0<1$) and $\zeta(1)=1$. For simplicity, we still denote $\tilde{F}$ by $F$.

By Galerkin method, it is not difficult to prove that above equation has a solution in $L^{\infty}([0,T],H^1(\mathbb{T},\mathfrak{g}))\bigcap W^{1,2}_2([0,T]\times\mathbb{T},\mathfrak{g})$ provided $\alpha>0$. The proof is similar with  that in the previous section. For the completeness, we provide a sketch of proof.

\subsection{Galerkin Approximatin of (\ref{lie:1})} Let $\lambda_i(i=1,2,...)$ be the eigenvalues of the operator $-f\Delta-\nabla f\cdot\nabla$ on the domain $H^2(\mathbb{T})$ and  $\{\omega_i: i=1, 2, \cdots\}$ is an orthonormal basis consisting of the eigenfunctions corresponding to $\lambda_i$. That is to say, for every $i\geq 1$,
\[-f\Delta\omega^i-\nabla f\cdot\nabla\omega^i=\lambda_i\omega^i.\]
The details of the eigenvalues of $-f\Delta-\nabla f\cdot\nabla$ can be found in chapter 2.4 of \cite{WMX}.

For the sake of convenience, we still denote
$$\mathfrak{J}(u)=\frac{u}{max\{|u|,1\}}.$$
According to Galerkin approximation, let
\begin{eqnarray*}
u^N(x,t):=\sum\limits_{i=1}^{N}\beta_i^N(t)\omega^i(x).
\end{eqnarray*}
Here $\{\beta_i^N(t)\}$ are unknown functions which take values in $\mathfrak{g}$ and assumed to satisfy the following ordinary differential equation:
\begin{equation}\label{2.1}
\left\{
\begin{array}{llll}
\aligned
&\alpha_0\frac{d\beta_i^N}{dt}+\alpha\sum\limits_{k=1}^N\int_{\mathbb{T}}\left[\mathfrak{J}(u^N),\,\frac{d\beta^N_k}{dt}\right]\omega^k\omega^i\,d\mathbb{T}
-\varepsilon\int_{\mathbb{T}}\sum\limits_{k=1}^{N}\beta_k^N(-\lambda_k\omega^k)\omega^i\,d\mathbb{T}\\
=&\int_{\mathbb{T}}\left[\mathfrak{J}(u^N),\,\sum\limits_{k=1}^{N}\beta_k^N(-\lambda_k\omega^k)\right]
\omega^i\,d\mathbb{T}+\int_{\mathbb{T}}F\left(t,x,\mathfrak{J}(u^N)\right)\omega^i\,d\mathbb{T},\\
&\beta_j^N(0)=\int_{\mathbb{T}}u_0\cdot\omega^j\,d\mathbb{T}.
\endaligned
\end{array}
\right.
\end{equation}
Since $F$ is $C^1$-smooth, the right side of (\ref{2.1}) is locally Lipschtiz continuous and there is a $\tau>0$ so that the solution of (\ref{2.1}) exists in $[0,\tau]$. Then, we get:
\begin{equation}\label{2.3}
\left\{
\begin{array}{llll}
  \aligned
  &\alpha_0\int_{\mathbb{T}}u^{N}_t\cdot\omega^i\,d\mathbb{T}+\alpha\int_{\mathbb{T}}[\mathfrak{J}(u^N),u^N_t]\omega^i d\mathbb{T}\\
  =&\,\varepsilon\int_{\mathbb{T}}(f\triangle u^N+\nabla f\cdot\nabla u^N)\omega^i\,d\mathbb{T}+\int_{\mathbb{T}}F(x,t,\mathfrak{J}(u^N))\omega^i\,d\mathbb{T}\\
  &+\int_{\mathbb{T}}[\mathfrak{J}(u^N),f\triangle u^N+\nabla f\cdot\nabla u^N]\omega^i\,d\mathbb{T},\\
  &u^N(0,\cdot)=u_0^N:=\sum\limits_{i=1}^N(\int_{\mathbb{T}}u_0\cdot\omega^i\,d\mathbb{T})\omega^i.
  \endaligned
\end{array}
\right.
\end{equation}
Multiplying the two sides of (\ref{2.3}) by $\beta_i^N(t)$ and summing $i$ from 1 to $N$ and integrating by parts, we get:
\[
\frac{\alpha_0}{2}\frac{d}{dt}\int_{\mathbb{T}}|u^N|^2\,d\mathbb{T}+\varepsilon\int_{\mathbb{T}}f|\nabla u^N|^2\,d\mathbb{T}=0,
\]
it follows
\begin{equation}\label{2.4}
\aligned
&\alpha_0\int_{\mathbb{T}}|u^N|^2\,d\mathbb{T}+2\varepsilon\int_0^t\int_{\mathbb{T}}f|\nabla u^N|^2\,d\mathbb{T}dt\\
=&\alpha_0\int_{\mathbb{T}}|u_0^N|^2\,d\mathbb{T}\leq \alpha_0\sum\limits_{i=1}^{\infty}\left(\int_{\mathbb{T}}u_0\cdot\omega^i\,d\mathbb{T}\right)^2
=\alpha_0\int_{\mathbb{T}}|u_0|^2\,d\mathbb{T}=\alpha_0\mbox{vol}(\mathbb{T}),
\endaligned
\end{equation}
where $\mbox{vol}(\mathbb{T})$ is the volume of $\mathbb{T}$.

Since
\[
\aligned
&\int_{\mathbb{T}}|u^N|^2\,d\mathbb{T}&=&\int_{\mathbb{T}}(\sum\limits_{i=1}^N\beta_i^N(t)\omega^i)(\sum\limits_{j=1}^N\beta_j^N(t)\omega^j)\,d\mathbb{T}\\
&&=&\sum\limits_{i=1}^N\sum\limits_{j=1}^N\beta_i^N(t)\beta_j^N(t)\cdot\delta_{ij}=\sum\limits_{i=1}^N|\beta_i^N(t)|^2,
\endaligned
\]
then, for any $T>0$ and any $i$, $\beta_i^N(t)$ can be extended to $[0,T]$. That is to say, $u^N$ can be extended to $[0,T]$.

Multiplying the two sides of (\ref{2.3}) by $-\lambda_i\beta_i^N$ and summing i from 1 to N, we get:
\[
\aligned
 &\int_{\mathbb{T}}\langle \alpha_0u^{N}_t, (f\triangle u^N+\nabla f\cdot\nabla u^N)\rangle\,d\mathbb{T}+\alpha\int_{\mathbb{T}}\langle\left[\mathfrak{J}(u^N),\, u^N_t\right], (f\Delta u^N+\nabla f\cdot\nabla u^N) \rangle d\mathbb{T}\\
  =&\varepsilon\int_{\mathbb{T}}|f\Delta u^N+\nabla f\cdot\nabla u^N|^2\,d\mathbb{T}
  +\int_{\mathbb{T}}\langle F\left(x,t,\mathfrak{J}(u^N)\right), (f\Delta u^N+\nabla f\cdot\nabla u^N)\rangle \,d\mathbb{T},
\endaligned
\]
this leads to
\begin{equation}\label{2.5}
\aligned
 &\frac{\alpha_0}{2}\frac{d}{dt}\int_{\mathbb{T}}f|\nabla u^N|^2\,d\mathbb{T}+\varepsilon\int_{\mathbb{T}}|f\Delta u^N+\nabla f\cdot\nabla u^N|^2\,d\mathbb{T}\\
=&\alpha\int_{\mathbb{T}}\langle\left[\mathfrak{J}(u^N),\, u^N_t\right], (f\Delta u^N+\nabla f\cdot\nabla u^N)\rangle d\mathbb{T}\\
  &-\int_{\mathbb{T}}\langle F\left(x,t,\mathfrak{J}(u^N)\right),(f\triangle u^N+\nabla f\cdot\nabla u^N)\rangle\,d\mathbb{T}.
\endaligned
\end{equation}
Multiplying the two sides of (\ref{2.3}) by $\frac{d\beta^N_i}{dt}$ and summing $i$ from $1$ to $N$ and integrating by parts, we get:
\begin{equation}\label{2.6}
\aligned
 &\alpha_0\int_{\mathbb{T}}|u^N_t|^2\,d\mathbb{T}+\frac{\varepsilon}{2}\frac{d}{dt}\int_{\mathbb{T}}f|\nabla u^N|^2\,d\mathbb{T}\\
=&-\int_{\mathbb{T}}\langle\left[\mathfrak{J}(u^N),\, u^N_t\right],\,(f\triangle u^N+\nabla f\cdot\nabla u^N) \rangle d\mathbb{T}+\int_{\mathbb{T}}F\left(x,t,\mathfrak{J}(u^N)\right)u_t^N\,d\mathbb{T}.
\endaligned
\end{equation}
Now, multiplying the two sides of (\ref{2.6}) by $\alpha$ and then adding the two sides of (\ref{2.5}) to the two sides of (\ref{2.6}) and integrating the obtained identity by parts, we deduce
\begin{equation}\label{vanish}
\aligned
&\frac{\alpha_0+\alpha\varepsilon}{2}\frac{d}{dt}\int_{\mathbb{T}}f|\nabla u^N|^2\,d\mathbb{T}+\varepsilon\int_{\mathbb{T}}|f\triangle u^N+\nabla f\cdot\nabla u^N|^2\,d\mathbb{T}+\alpha_0\alpha\int_{\mathbb{T}}|u_t^N|^2\,d\mathbb{T}\\
=&\alpha\int_{\mathbb{T}}\langle F\left(x,t,\mathfrak{J}(u^N)\right),\,u_t^N\rangle\,d\mathbb{T}+\int_{\mathbb{T}}\nabla\left( F\left(x,t,\mathfrak{J}(u^N)\right)\right)\cdot\nabla u^N\cdot f\,d\mathbb{T}.
\endaligned
\end{equation}
Note that
\[
\aligned
&|\nabla(F(x,t,\mathfrak{J}(u^N)))|\leq|(\nabla_x F)(x,t,\mathfrak{J}(u^N))| +\left| \frac{\partial F}{\partial z}\left(x,t,\mathfrak{J}(u^N)\right)\cdot\nabla\left(\mathfrak{J}(u^N)\right)\right|,
\endaligned
\]

\[
\aligned
\nabla\left(\mathfrak{J}(u^N)\right)=\left(\nabla\left(\frac{u_1^N}{max\{|u^N|,1\}}\right),\, ...,\, \nabla\left(\frac{u_m^N}{max\{|u^N|,1\}}\right)\right),
\endaligned
\]
and for $k=1,2,...,m$
\[
\aligned
\nabla\Bigg(\frac{u_k^N}{max\{|u^N|,1\}}\Bigg)
=\frac{\nabla u_k^N}{max\{|u^N|,1\}}-\chi_{\{|u^N|>1\}}\frac{u_k^N(\sum\limits_{i=1}^m\nabla u_i^N\cdot u_i^N)}{|u^N|^3}
\endaligned
\]
where
\[
\chi_{\{|u^N|>1\}}(x)=\Big\{\begin{array}{cc}0,& \s\s |u^N(x)|\leq 1,\\
1,&\s\s|u^N(x)|>1.
\end{array}
\]
So, there exists a $C_1$ such that
\[
|\nabla(F(x,t,\mathfrak{J}(u^N)))|\leq C_1|\nabla u^N|+C_1.
\]
Obviously, there holds true
\[
|F(x,t,\mathfrak{J}(u^N))|\leq C_2.
\]
In view of the above estimates we have
\[
\aligned
&\frac{\alpha_0+\alpha\varepsilon}{2}\frac{d}{dt}\int_{\mathbb{T}}f|\nabla u^N|^2\,d\mathbb{T}+\varepsilon\int_{\mathbb{T}}|f\triangle u^N+\nabla f\cdot\nabla u^N|^2\,d\mathbb{T}+\alpha_0\alpha\int_{\mathbb{T}}|u_t^N|^2\,d\mathbb{T}\\
\leq&\alpha C_2\int_{\mathbb{T}}|u_t^N|\,d\mathbb{T}+\int_{\mathbb{T}}(C_1+C_1|\nabla u^N|)\cdot|\nabla u^N|\cdot f\,d\mathbb{T}\\
\leq&\frac{\alpha\alpha_0}{2}\int_{\mathbb{T}}|u_t^N|^2\,d\mathbb{T}+\int_{\mathbb{T}}\frac{\alpha C_2^2}{2\alpha_0}\,d\mathbb{T}+C_1\int_{\mathbb{T}}(|\nabla u^N|+|\nabla u^N|^2)f\,d\mathbb{T}\\
\leq&\frac{\alpha\alpha_0}{2}\int_{\mathbb{T}}|u_t^N|^2\,d\mathbb{T}+\frac{\alpha C_2^2}{2\alpha_0}\mbox{vol}(\mathbb{T})+C_1\int_{\mathbb{T}}(\frac{1+|\nabla u^N|^2}{2}+|\nabla u^N|^2)f\,d\mathbb{T}\\
\leq&\frac{\alpha\alpha_0}{2}\int_{\mathbb{T}}|u_t^N|^2\,d\mathbb{T}+\frac{\alpha C_2^2}{2\alpha_0}\mbox{vol}(\mathbb{T})+\frac{3}{2}C_1\int_{\mathbb{T}}(1+|\nabla u^N|^2)f\,d\mathbb{T}\\
\leq&\frac{\alpha\alpha_0}{2}\int_{\mathbb{T}}|u_t^N|^2\,d\mathbb{T}+(\frac{\alpha C_2^2}{2\alpha_0}+\frac{3C_1}{2}||f||_{\infty})\mbox{vol}(\mathbb{T})+\frac{3}{2}C_1\int_{\mathbb{T}}f|\nabla u^N|^2\,d\mathbb{T}.
\endaligned
\]
After rearranging the terms in the above inequality we have
\begin{equation}\label{2.7}
\aligned
&\frac{\alpha_0+\alpha\varepsilon}{2}\int_{\mathbb{T}}f|\nabla u^N|^2\,d\mathbb{T}+\varepsilon\int_0^t\int_{\mathbb{T}}|\Delta_f u^N|^2\,d\mathbb{T}dt +\frac{\alpha_0\alpha}{2}\int_0^t\int_{\mathbb{T}}|u_t^N|^2\,d\mathbb{T}dt\\
\leq&(\frac{\alpha C_2^2}{2\alpha_0}+\frac{3C_1}{2}||f||_{\infty})\cdot\mbox{vol}(\mathbb{T})t+\frac{3}{2}C_1\int_0^t\int_{\mathbb{T}}f|\nabla u^N|^2\,d\mathbb{T}dt\\
&+\frac{\alpha_0+\alpha\varepsilon}{2}\int_{\mathbb{T}}f|\nabla u^N_0|^2\,d\mathbb{T}.
\endaligned
\end{equation}
On the other hand, we have
\[
\aligned
&\int_{\mathbb{T}}f\cdot|\nabla u_0|^2\,d\mathbb{T}=\int_{\mathbb{T}}(-f\triangle u_0-\nabla f\cdot\nabla u_0)\cdot u_0\,d\mathbb{T}\\
=&\int_{\mathbb{T}}\sum\limits_{i=1}^{\infty}\left\{\int_{\mathbb{T}}u_0\cdot(-f\triangle\omega^i-\nabla f\cdot\nabla\omega^i)\,d\mathbb{T}\right\}\omega^i\cdot u_0\,d\mathbb{T}=\sum\limits_{i=1}^{\infty}\lambda_i\left|\int_{\mathbb{T}}u_0\cdot\omega^i\,d\mathbb{T}\right|^2
\endaligned
\]
and
\[
\aligned
&\int_{\mathbb{T}}f|\nabla u^N_0|^2\,d\mathbb{T}&=&\int_{\mathbb{T}}(-f\Delta u_0^N-\nabla f\cdot\nabla u_0^N)\cdot u_0^N\,d\mathbb{T}=\sum\limits_{i=1}^{N}\lambda_i\left|\int_{\mathbb{T}}u_0\cdot\omega^i\,d\mathbb{T}\right|^2.
\endaligned
\]

Let $$E(t)=\int_0^t\int_{\mathbb{T}}f|\nabla u^N|^2\,d\mathbb{T}dt,$$
we infer from (\ref{2.7}) that
\[
\frac{\alpha_0+\alpha\varepsilon}{2}\frac{d}{dt}E\leq(\frac{\alpha C_2^2}{2\alpha_0}+\frac{3C_1}{2}||f||_{\infty})\mbox{vol}(\mathbb{T})t+\frac{3}{2}C_1E+\frac{\alpha_0+\alpha\varepsilon}{2}\int_{\mathbb{T}}f\cdot|\nabla u_0|^2\,d\mathbb{T}.
\]
By Gronwall inequality, we obtain:
\[
E(t)\leq C_8(T).
\]
Here $C_8(T)$ is independent of $\varepsilon$. So, there holds true
\begin{eqnarray}\label{ZL}
&&\frac{\alpha_0+\alpha\varepsilon}{2}\int_{\mathbb{T}}f|\nabla u^N(t)|^2\,d\mathbb{T}+\varepsilon\int_0^t\int_{\mathbb{T}}|\Delta_f u^N|^2 \,d\mathbb{T}dt +\frac{\alpha\alpha_0}{2}\int_0^t\int_{\mathbb{T}}|u_t^N|^2\,d\mathbb{T}dt\\
&\leq&(\frac{\alpha C_2^2}{2\alpha_0}+\frac{3C_1}{2}||f||_{\infty})\mbox{vol}(\mathbb{T})t+\frac{3}{2}C_1\cdot C_8(T)+\frac{\alpha_0+\alpha\varepsilon}{2}\int_{\mathbb{T}}f\cdot|\nabla u_0|^2\,d\mathbb{T}.\nonumber
\end{eqnarray}
So, we have also
\begin{equation}\label{alpha}
\alpha\alpha_0\int_0^t\int_{\mathbb{T}}|u_t^N|^2\,d\mathbb{T}dt\leq C(\frac{1}{\alpha_0}, \alpha).
\end{equation}

On the other hand, as $0<\tilde{m}\leq f(x)\leq ||f||_{\infty}<\infty$, the above inequality implies
\[
\aligned
&\frac{\alpha_0+\alpha\varepsilon}{2}\tilde{m}\int_{\mathbb{T}}|\nabla u^N(t)|^2\,d\mathbb{T}
\leq\frac{\alpha_0+\alpha\varepsilon}{2}\int_{\mathbb{T}}f|\nabla u^N(t)|^2\,d\mathbb{T}\\
\leq&(\frac{\alpha C_2^2}{2\alpha_0}+\frac{3C_1}{2}||f||_{\infty})\mbox{vol}(\mathbb{T})t+\frac{3}{2}C_1\cdot C_8(T)+\frac{\alpha_0+\alpha\varepsilon}{2}\int_{\mathbb{T}}f\cdot|\nabla u_0|^2\,d\mathbb{T}.
\endaligned
\]
It follows that there exists a constant $C_{10}(T)$ which is independent of $\varepsilon$ such that
\[
\int_{\mathbb{T}}|\nabla u^N|^2\,d\mathbb{T}\leq C_{10}(T)
\]
and
\begin{eqnarray*}
\begin{array}{lll}
\dint_0^t\dint_{\mathbb{T}}|\triangle u^N|^2\,d\mathbb{T}dt &\leq& 2\dint_0^t\dint_{\mathbb{T}}\bigg|\Delta u^N+\frac{\nabla f}{f}\nabla u^N\bigg|^2\,d\mathbb{T}dt+2\dint_0^t\dint_{\mathbb{T}}\left|-\frac{\nabla f}{f}\nabla u^N\right|^2\,d\mathbb{T}dt\\
&\leq& \displaystyle{\frac{2}{\tilde{m}^2}}\dint_0^t\dint_{\mathbb{T}}|\Delta_f u^N|^2\,d\mathbb{T}dt+2\left\|\frac{\nabla f}{f}\right\|^2_{\infty}\int_0^t\int_{\mathbb{T}}|\nabla u^N|^2\,d\mathbb{T}dt\\
&\leq& \displaystyle{\frac{1+\alpha\varepsilon}{\varepsilon \tilde{m}^2}}\dint_{\mathbb{T}}f\cdot|\nabla u_0|^2\,d\mathbb{T}+2\left\|\frac{\nabla f}{f}\right\|^2_{\infty}\frac{1}{\tilde{m}}\int_0^t\int_{\mathbb{T}}f|\nabla u^N|^2\,d\mathbb{T}dt\\
& &+\left( \displaystyle{\frac{\alpha C_2^2}{\varepsilon\tilde{m}^2}} +  \displaystyle{\frac{3C_1}{\varepsilon \tilde{m}^2}}||f||_{\infty}\right)\mbox{vol}(\mathbb{T})t+ \displaystyle{\frac{3}{\varepsilon \tilde{m}^2}}C_1\cdot C_8(T)\\
&\leq&\left( \displaystyle{\frac{\alpha C_2^2}{\varepsilon\tilde{m}^2}} +  \displaystyle{\frac{3C_1}{\varepsilon \tilde{m}^2}}||f||_{\infty}\right)\mbox{vol}(\mathbb{T})t+ \displaystyle{\frac{3}{\varepsilon \tilde{m}^2}}C_1\cdot C_8(T)\\
& &+ \displaystyle{\frac{1+\alpha\varepsilon}{\varepsilon \tilde{m}^2}}\dint_{\mathbb{T}}f\cdot|\nabla u_0|^2\,d\mathbb{T}+2\left\|\frac{\nabla f}{f}\right\|^2_{\infty}\frac{1}{\tilde{m}}C_8(T).
\end{array}
\end{eqnarray*}
So we get

\begin{lem} Assume that $f$ and $F$ satisfy respectively the same conditions as in Theorem \ref{thm1}. Then, the approximate solution sequence $\{u^N\}$ to (\ref{2.1}) satisfies

$\bullet$ $\{u^N\}$ is a bounded sequence in $L^{\infty}([0,T],H^1(\mathbb{T},\mathfrak{g}))$;

\medskip
$\bullet$ $\{u^N_t\}$ is a bounded sequence in $L^{2}([0,T],L^2(\mathbb{T},\mathfrak{g}))$;

\medskip
$\bullet$ $\{\Delta u^N\}$ is a bounded sequence in $L^{2}([0,T],L^2(\mathbb{T},\mathfrak{g}))$;

\medskip
$\bullet$ $\{\nabla u^N\}$ is a bounded sequence in $L^{\infty}([0,T],L^2(\mathbb{T},T\mathbb{T}\otimes\mathfrak{g}))$,  where $T\mathbb{T}$ is the tangent bundle of $\mathbb{T}$.
\end{lem}

\subsection{Derivation of auxilliary limit equation}
By the property of weak limits and Aubin-Lions Lemma, there exists a $v^{\varepsilon}\in W^{2,1}_2(\mathbb{T}\times[0,T],\mathfrak{g})$ and a subsequence of $\{u^N\}$ which is also denoted by $\{u^N\}$ such that

\medskip
$\bullet$ $u^N\rightharpoonup v^{\varepsilon}\s \mbox{weakly* in}\,\,L^{\infty}([0,T],H^1(\mathbb{T},\mathfrak{g}))$;

\medskip
$\bullet$ $u^N\rightarrow v^{\varepsilon}\s\mbox{strongly in}\,\,L^{\infty}([0,T],L^2(\mathbb{T},\mathfrak{g}))$;

\medskip
$\bullet$ $u^N\rightarrow v^{\varepsilon}\,\,\,\,a.e.\,\,\mathbb{T}\times[0,T]$;

\medskip
$\bullet$ $u^N_t\rightharpoonup v^{\varepsilon}_t\,\,\,\,\mbox{weakly in}\,\,L^{2}([0,T],L^2(\mathbb{T},\mathfrak{g}))$;

\medskip
$\bullet$ $\Delta u^N\rightharpoonup \Delta v^{\varepsilon}\,\,\,\,\mbox{weakly in}\,\,L^2([0,T],L^2(\mathbb{T},\mathfrak{g}))$;

\medskip
$\bullet$ $\nabla u^N\rightharpoonup \nabla v^{\varepsilon}\,\,\,\,\mbox{weakly* in}\,\,L^{\infty}([0,T],L^2(\mathbb{T},T\mathbb{T}\otimes\mathfrak{g}))$.

\medskip
\noindent Since
\[
||u^N||_{L^{\infty}([0,T],H^1(\mathbb{T},\mathfrak{g}))}\leq C_{12}
\]
and $u^N\rightharpoonup v^{\varepsilon}\,\,\mbox{weakly* in}\,\,L^{\infty}([0,T],H^1(\mathbb{T},\mathfrak{g}))$, we have
\begin{eqnarray}\label{ZL2}
||v^{\varepsilon}||_{L^{\infty}([0,T],H^1(\mathbb{T},\mathfrak{g}))}\leq C_{12}.
\end{eqnarray}
By the same method as we prove $(\ref{ZL2})$, from (\ref{ZL}) it follows that
\begin{equation}\label{vanish1}
\alpha\alpha_0\int_0^T\int_{\mathbb{T}}|v_t^{\varepsilon}|^2\,d\mathbb{T}\,dt\leq C_{13},
\end{equation}
where $C_{13}\equiv  C_{13}(\frac{1}{\alpha_0}, T)$ depends on $\frac{1}{\alpha_0}$ and $T$. So, it is easy to see that
\[
\left[\mathfrak{J}(u^N),u^N_t\right]\rightharpoonup\left[\mathfrak{J}(v^{\varepsilon}),v^{\varepsilon}_t\right]
\]
weakly in $L^{2}([0,T],L^2(\mathbb{T},\mathfrak{g}))$,
\[
\left[\mathfrak{J}(u^N),\Delta u^N\right]\rightharpoonup\left[\mathfrak{J}(v^{\varepsilon}),\Delta v^{\varepsilon}\right]
\]
weakly in $L^{2}([0,T],L^2(\mathbb{T},\mathfrak{g}))$, and
\[
\left[\mathfrak{J}(u^N),\nabla f\cdot\nabla u^N\right]\rightharpoonup\left[\mathfrak{J}(v^{\varepsilon}),\nabla f\cdot\nabla v^{\varepsilon}\right]
\]
weakly$*$ in $L^{\infty}([0,T],L^2(\mathbb{T},\mathfrak{g}))$.

Fix $r\in \mathbb{Z}^+$ and take any $N\geq r$. First, we multiply two sides of (\ref{2.3}) by $\eta^i(t)$ which belongs to $C^{\infty}([0,T],\mathfrak{g})$ and sum $i$ from 1 to r,  then integrate the obtained identity on $[0,T]$ to derive
\[
\aligned
&\alpha_0\int_0^T\int_{\mathbb{T}}u^{N}_t\cdot\Phi^r\,d\mathbb{T}dt+\alpha\int_0^T\int_{\mathbb{T}}\left[\mathfrak{J}(u^N),\,u^N_t\right]\Phi^r d\mathbb{T}dt\\
  =&\varepsilon\int_0^T\int_{\mathbb{T}}(f\triangle u^N+\nabla f\cdot\nabla u^N)\Phi^r\,d\mathbb{T}dt+\int_0^T\int_{\mathbb{T}}\left[\mathfrak{J}(u^N),f\triangle u^N+\nabla f\cdot\nabla u^N\right]\Phi^r\,d\mathbb{T}dt\\
  &+\int_0^T\int_{\mathbb{T}}F\left(x,t,\mathfrak{J}(u^N)\right)\Phi^r\,d\mathbb{T}dt.
\endaligned
\]
Here
\[
\Phi^r(x,t)=\sum\limits_{i=1}^r\omega^i(x)\eta^i(t).
\]

Letting $N$ tends to $\infty$ in the above identity, we get:
\[
\aligned
&\alpha_0\int_0^T\int_{\mathbb{T}}v_t^{\varepsilon}\cdot\Phi^r\,d\mathbb{T}dt+\alpha\int_0^T\int_{\mathbb{T}}\left[\mathfrak{J}(u^N),\,v_t^{\varepsilon}\right]\Phi^r\,d\mathbb{T}dt\\
=&\int_0^T\int_{\mathbb{T}}[\mathfrak{J}(u^N),\,f\Delta v^{\varepsilon}+\nabla f\cdot\nabla v^{\varepsilon}]\Phi^r\,d\mathbb{T}dt +\int_0^T\int_{\mathbb{T}}F\left(x,t,\mathfrak{J}(u^N)\right)\Phi^r\,d\mathbb{T}dt\\
&+\varepsilon\int_0^T\int_{\mathbb{T}}(f\triangle v^{\varepsilon}+\nabla f\cdot\nabla v^{\varepsilon}
)\Phi^r\,d\mathbb{T}dt.
\endaligned
\]
Since the functions, which are of type $\Phi^r(x,t)$, are dense in $L^2([0,T],L^2(\mathbb{T},\mathfrak{g}))$, we conclude that, in the sense of distribution, there holds true
\begin{eqnarray}\label{Key}
\left\{
\begin{array}{llll}
  \aligned
 & \alpha_0v^{\varepsilon}_t+\alpha[\mathfrak{J}(v^{\varepsilon}),v^{\varepsilon}_t]=&\varepsilon(f\triangle v^{\varepsilon}+\nabla f\cdot\nabla v^{\varepsilon})+F(x,t,\mathfrak{J}(v^{\varepsilon}))\\
 &  &+[\mathfrak{J}(v^{\varepsilon}),\, f\triangle v^{\varepsilon}+\nabla f\cdot\nabla v^{\varepsilon}],\s\s\s\s\s\\
 &v^{\varepsilon}(0,\cdot)=u_0.&
\endaligned
\end{array}
\right.
\end{eqnarray}

Now, respectively we still choose $$v^{\varepsilon}-v^{\varepsilon}\frac{min\{1,|v^{\varepsilon}|\}}{|v^{\varepsilon}|}$$
and
$$\frac{v^{\varepsilon}(max\{|v^{\varepsilon}|,1\}-1)}{|v^{\varepsilon}|(|v^{\varepsilon}|-1+\delta)}$$ as test functions of the above equation (\ref{Key}). By almost the same calculation as in the above section we obtain
\[
\frac{d}{dt}\int_{|v^{\varepsilon}|>1}|v^{\varepsilon}|^2\Big(1-\frac{1}{|v^{\varepsilon}|}\Big)\,d\mathbb{T}\leq 0.
\]
This means that the following function
\[
q(t):=\int_{|v^{\varepsilon}|>1}|v^{\varepsilon}|^2\Big(1-\frac{1}{|v^{\varepsilon}|}\Big)\,d\mathbb{T}
\]
is decreasing non-negative function. Noting $|v^{\varepsilon}(\cdot,0)|=|u_0|=1$, i.e. $q(0)=0$, we can see that $q(t)\equiv0$ for any $t>0$. Therefore, for any $t\in[0, T]$ we have
$$|v^{\varepsilon}(t)|\leq 1 \s\s\s a.\, e.\s\mathbb{T}.$$
Immediately it follows that $v^{\varepsilon}$ satisfies in the sense of distribution
\[
\alpha_0v^{\varepsilon}_t+\alpha[v^{\varepsilon},\,v^{\varepsilon}_t]=\varepsilon(f\triangle v^{\varepsilon}+\nabla f\cdot\nabla v^{\varepsilon})+[v^{\varepsilon},\,f\triangle v^{\varepsilon}+\nabla f\cdot\nabla v^{\varepsilon}]+F(x,t,v^{\varepsilon})
\]
with initial value $v^{\varepsilon}(\cdot, 0)=u_0$.

\subsection{The Proof of Theorem \ref{thm1}} We need to discuss the two cases:
\medskip

{\bf Case 1.  $\alpha_0>0$ and $\alpha>0$.}
Multiplying the two sides of the above equation by $v^{\varepsilon}$ and integrating it on $\mathbb{T}\times[0,t]$, we get:
\begin{equation}\label{2.8}
\alpha_0\int_{\mathbb{T}}(|v^{\varepsilon}(t)|^2-1)\,d\mathbb{T}+2\varepsilon\int_0^t\int_{\mathbb{T}}f\cdot|\nabla v^{\varepsilon}|^2\,d\mathbb{T}d\tau=0.
\end{equation}
Note that $\{v^{\varepsilon}\}$ is a bounded sequence in $L^{\infty}([0,T],H^1(\mathbb{T},\mathfrak{g}))$ and $\{v_t^{\varepsilon}\}$ is a bounded sequence in $L^2([0,T],L^2(\mathbb{T},\mathfrak{g}))$. So, by the property of weak limits and Aubin-Lions Lemma, there is a $u$ and a subsequence of $\{v^{\varepsilon}\}$ which is also denoted by $\{v^{\varepsilon}\}$ such that:

$\bullet$ $v^{\varepsilon}\rightharpoonup u$ weakly $*$ in $L^{\infty}([0,T],H^1(\mathbb{T},\mathfrak{g}))$;
\medskip

$\bullet$ $v^{\varepsilon}\rightarrow u$ strongly in $L^{\infty}([0,T],L^2(\mathbb{T},\mathfrak{g}))$;

\medskip
$\bullet$ $v^{\varepsilon}_t\rightharpoonup u_t$ weakly in $L^{2}([0,T],L^2(\mathbb{T},\mathfrak{g}))$.

\medskip
Letting $\varepsilon$ tends to 0 in (\ref{2.8}), we have that, for all $t\in[0,T]$, there holds true
\[
\int_{\mathbb{T}}(|u|^2-1)\,d\mathbb{T}=0.
\]
This leads to $$|u|=1\s\s\s \mbox{for}\s a.e.\s (x, t)\in\mathbb{T}\times[0, T],$$
since $|v^{\varepsilon}|\leq1$ implies that $|u|\leq1$.

From $v^{\varepsilon}\rightarrow u$ $a.e.$ on $\mathbb{T}\times[0,T]$ we can deduce easily that $F(x,t,v^{\varepsilon})\rightarrow F(x,t,u)$ $a.e.$ on $\mathbb{T}\times[0,T]$.

For any $\varphi\in C^{\infty}(\mathbb{T}\times[0,T],\mathfrak{g})$ we have
\[
\aligned
&\alpha_0\int_{\mathbb{T}}v^{\varepsilon}(T)\cdot\varphi(T)\,d\mathbb{T}-\alpha_0\int_{\mathbb{T}}u_0\cdot\varphi(0)\,d\mathbb{T}
+\alpha\int_0^T\int_{\mathbb{T}}[v^{\varepsilon},v^{\varepsilon}_t]\cdot\varphi\,d\mathbb{T}\,dt\\
=&\alpha_0\int_{0}^{T}\int_{\mathbb{T}}v^{\varepsilon}\varphi_t\,d\mathbb{T}\,dt-\int_0^T\int_{\mathbb{T}}[v^{\varepsilon},f\nabla v^{\varepsilon}]\cdot\nabla \varphi\,d\mathbb{T}\,dt+\int_0^T\int_{\mathbb{T}}F(x,t,v^{\varepsilon})\cdot\varphi\,d\mathbb{T}\,dt\\
&-\varepsilon\int_0^T\int_{\mathbb{T}}f\nabla v^{\varepsilon}\cdot\nabla \varphi\,d\mathbb{T}\,dt.
\endaligned
\]
Hence, it follows
\begin{equation}\label{Fin}
\aligned
&\alpha_0\int_{\mathbb{T}}u(T)\varphi(T)\,d\mathbb{T}-\alpha_0\int_{\mathbb{T}}u_0\varphi(0)\,d\mathbb{T}
+\alpha\int_{0}^{T}\int_{\mathbb{T}}[u,u_t]\varphi\,d\mathbb{T}dt\\
=&\alpha_0\int_{0}^{T}\int_{\mathbb{T}}u\varphi_t\,d\mathbb{T}dt-\int_{0}^{T}\int_{\mathbb{T}}[u,f\nabla u]\nabla\varphi\,d\mathbb{T}dt+\int_{0}^{T}\int_{\mathbb{T}}F(x,t,u)\varphi\,d\mathbb{T}dt,
\endaligned
\end{equation}
since the controlled convergence theorem tells us that, as $\varepsilon\rightarrow 0$,
\[
\int_0^T\int_{\mathbb{T}}F(x,t,v^{\varepsilon})\cdot\varphi\,d\mathbb{T}dt\rightarrow\int_0^T\int_{\mathbb{T}}F(x,t,u)\cdot\varphi\,d\mathbb{T}dt.
\]

\medskip
{\bf Case 2. $\alpha_0>0$ and $\alpha=0$.} We have known that, as $\alpha>0$, $\{v^{\varepsilon}\}$ is a bounded sequence in the Sobolev space $L^{\infty}([0,T],H^1(\mathbb{T},\mathfrak{g}))\cap W^{1,1}_2(\mathbb{T}\times[0,T], \mathfrak{g})$. By the same argument as in the above, letting $\varepsilon$ in (\ref{2.8}) tends to $0$ and denoting the limit of $v^\varepsilon$ by $u^\alpha$ we conclude that, for any $\varphi\in C^{\infty}([0,T]\times\mathbb{T},\mathfrak{g})$,
\begin{equation}\label{Fin'}
\aligned
&\alpha_0\int_{\mathbb{T}}u^{\alpha}(T)\varphi(T)\,d\mathbb{T}-\alpha_0\int_{\mathbb{T}}u_0\varphi(0)\,d\mathbb{T}
+\alpha\int_{0}^{T}\int_{\mathbb{T}}[u^{\alpha},u^{\alpha}_t]\varphi\,d\mathbb{T}dt\\
=&\alpha_0\int_{0}^{T}\int_{\mathbb{T}}u^{\alpha}\varphi_t\,d\mathbb{T}dt-\int_{0}^{T}\int_{\mathbb{T}}[u^{\alpha},f\nabla u^{\alpha}]\nabla\varphi\,d\mathbb{T}dt+\int_{0}^{T}\int_{\mathbb{T}}F(x,t,u^{\alpha})\varphi\,d\mathbb{T}dt.
\endaligned
\end{equation}
Obviously, we have $u^\alpha$ is uniformly bounded in $L^{\infty}([0,T],H^1(\mathbb{T},\mathfrak{g}))$. Therefore, there exists a $u\in L^{\infty}([0,T],H^1(\mathbb{T},\mathfrak{g}))$ and a subsequence of $u^{\alpha}$, which is still denoted by $u^{\alpha}$, such that $u^\alpha\to u$ weakly $\ast$ in $L^{\infty}([0,T],H^1(\mathbb{T},\mathfrak{g}))$ and $u^\alpha\to u$ a.e. $\mathbb{T}\times[0,T]$.

From (\ref{vanish1}) we know that, as $\alpha>0$,
\begin{equation*}
\alpha\alpha_0\int_0^t\int_{\mathbb{T}}|u^{\alpha}_t|^2\,d\mathbb{T}dt\leq C(\frac{1}{\alpha_0}, \alpha, T).
\end{equation*}
It follows that, as $\alpha \rightarrow 0$, there holds true
$$\alpha\int_{0}^{T}\int_{\mathbb{T}}[u^{\alpha},u^{\alpha}_t]\varphi\,d\mathbb{T}dt\longrightarrow 0.$$
Letting $\alpha \rightarrow 0$ in (\ref{Fin'}) we get
\begin{equation*}
\aligned
&\alpha_0\int_{\mathbb{T}}u(T)\varphi(T)\,d\mathbb{T}-\alpha_0\int_{\mathbb{T}}u_0\varphi(0)\,d\mathbb{T}\\
=&\alpha_0\int_{0}^{T}\int_{\mathbb{T}}u\varphi_t\,d\mathbb{T}dt-\int_{0}^{T}\int_{\mathbb{T}}[u,f\nabla u]\nabla\varphi\,d\mathbb{T}dt+\int_{0}^{T}\int_{\mathbb{T}}F(x,t,u)\varphi\,d\mathbb{T}dt.
\endaligned
\end{equation*}

Up to now, we have proved the following
\begin{thm}
Let $(\mathbb{T}, h)$ be an $n$-dimensional closed manifolds equipped with a metric $h$ and $\mathfrak{g}$ be a $m$-dimensional compact Lie algebra. Assume that $\alpha_0>0$, $\alpha \geq 0$, $F(x, t, z): \mathbb{T}\times\mathbb{R}^+\times S_{\mathfrak{g}}(1)\to \mathfrak{g}$ is $C^1$-smooth and $f\in C^1(\mathbb{T})$ with $\min_{x\in\mathbb{T}}f(x)>0$. Then, in the case $\alpha_0>0$ and $\alpha>0$ (\ref{Popu}) admits a global weak solution $u\in W^{1,1}_2(\mathbb{T}\times[0, T], S_{\mathfrak{g}}(1))$ for any $T>0$; in the case $\alpha_0>0$ and $\alpha=0$ (\ref{Popu}) admits a alobal weak solution $u\in L^\infty_{loc}(\mathbb{R}^+,H^1(\mathbb{T}, S_{\mathfrak{g}}(1)))$; provided the initial value map $u_0$ belongs to $H^{1}(\mathbb{T}, S_{\mathfrak{g}}(1))$.
\end{thm}

Since Theorem \ref{thm1} is a special case of the above theorem, thus we also complete the proof of Theorem \ref{thm1}.

\begin{rem}
From the above arguments, we can see easily that the Cauchy problem for (\ref{Popu}) admits a global weak solution if the coupling function $f=f(x, t)>0$ defined on $\mathbb{T}\times [0, \infty)$ depends on time variable and is $C^1$-smooth with respect to $x$ and $t$. In fact, such a system is also of physical background when $\mathfrak{g}=\mathbb{R}^3$, for more details we refer to \cite{JW} and references therein.
\end{rem}

%\vspace{1cm}

{}

\vspace{0.5cm}

Zonglin Jia

{\small\it Institute of Applied Physics and Computational Mathematics, China Academy of Engineering Physics, Beijing, 100088, P. R. China}

{\small\it Email: 756693084@qq.com}\\

Youde Wang

{\small \it College of Mathematics and Information Sciences, Guangzhou University.}
\medskip

{\small\it  Academy of Mathematics and Systems Science, Chinese Academy of Sciences, Beijing 100190,  P. R. China.}

{\small\it  Email: wyd@math.ac.cn}

\end{document}